\documentclass[12pt,oneside, reqno]{amsart}
\usepackage{amssymb,latexsym,cite}
\usepackage[usenames, dvipsnames]{color}
\definecolor{mygray}{gray}{0.6}

\usepackage{amssymb,amsfonts,amsthm, amsmath, breqn, color, float, mathtools, booktabs, hvfloat, url,cite}

\advance\textheight by \topskip   \numberwithin{equation}{section}

\newtheorem{theorem}{Theorem}[section]

\newtheorem{lemma}[theorem]{Lemma}

\newcommand{\beq}{\begin{eqnarray*}}
\newcommand{\feq}{\end{eqnarray*}}
\newcommand{\beqn}{\begin{eqnarray}}
\newcommand{\feqn}{\end{eqnarray}}

\def\obrace{\iftrue{\else}\fi}
\def\cbrace{\iffalse{\else}\fi}

\let\originalparagraph\paragraph
\renewcommand{\paragraph}[2][.]{\originalparagraph{#2#1}}

\newcommand{\calt}{{\mathcal T}}

\newcommand{\bes}{\begin{split}}
\newcommand{\fes}{\end{split}}

\newtheorem*{conj*}{Conjecture}
\makeatletter \@addtoreset{theorem}{section}\makeatother

\makeatletter \@addtoreset{theorem}{section}\makeatother

\makeatletter \@addtoreset{theorem}{section}\makeatother

\newtheorem*{theorem*}{Theorem}

\makeatletter
\newcommand{\leqnomode}{\tagsleft@true\let\veqno\@@leqno}
\newcommand{\reqnomode}{\tagsleft@false\let\veqno\@@eqno}
\makeatother
\usepackage{algorithm}
\usepackage{graphicx}
\usepackage[noend]{algpseudocode}
\makeatletter
\def\BState{\State\hskip-\ALG@thistlm}
\makeatother
\DeclarePairedDelimiter\ceil{\lceil}{\rceil}

\DeclareCaptionLabelFormat{noname}{#2}
\newlength\myindent
\def\B{\mathcal{B}}
\def\CP{\mathcal{CP}}
\def\CCP{\mathcal{CCP}}
\begin{document}
\title[Baryiamonds and Polyiamonds]{Enumeration of Various Animals on the triangular lattice}
\author{Toufik Mansour}
\address{Department of Mathematics, University of Haifa, 3498838 Haifa, Israel}
\email{tmansour@univ.haifa.ac.il }

\author{Reza Rastegar}
\address{Occidental Petroleum Corporation, Houston, TX 77046 and Departments of Mathematics and Engineering, University of Tulsa, OK 74104, USA - Adjunct Professor}
\email{reza\_rastegar2@oxy.com}


\begin{abstract}
In this paper, we consider various classes of polyiamonds that are animals residing on the triangular lattice. By careful analyses through certain layer-by-layer decompositions and cell pruning/growing arguments, we derive explicit forms for the generating functions of the number of nonempty translation-invariant baryiamonds (bargraphs in the triangular lattice), column-convex polyiamonds, and convex polyiamonds with respect to their perimeter. In particular, we show that the number of \\
(A) baryiamonds of perimeter $n$ is asymptotically
$$\frac{(\xi+1)^2\sqrt{\xi^4+\xi^3-2\xi+1}}{2\sqrt{\pi n^3}}\xi^{-n-2},$$
where $\xi$ is a root of a certain explicit polynomial of degree 5. \\
(B) column-convex polyiamonds of perimeter $n$ is asymptotic to
$$\frac{(17997809\sqrt{17}+3^3\cdot13\cdot175463)\sqrt{95\sqrt{17}-119}}{2^7\cdot43^2\cdot 89^2\sqrt{6\pi n^3}}\left(\frac{3+\sqrt{17}}{2}\right)^{n-1}.$$\\
(C) convex polyiamonds of perimeter $n$ is asymptotic to
$$\frac{1280}{441\sqrt{3\pi n^3}}3^n.$$
\end{abstract}
\maketitle

\section{Introduction}

An animal, living on a two dimensional lattice, is an edge-connected set of basic two-dimensional polygon-like cells on a two-dimension lattice, where the connectivity of two cells is defined by having a common edge. Among several classes of animals are squared polyominoes or simply {\it polyominoes}. These are animals that are made by edge-connected squares on the squared lattice and are well studied in statistical physics, combinatorics, and discrete geometry. They specifically play a role in enumeration of graphs  \cite{har1}, modeling the mechanics of macro-molecules \cite{Temperly} such as the collapse of branched polymers \cite{peard1}, percolation processes \cite{bro1}, and cell growth processes \cite{eden1, klar1, Read}, to name a few. In combinatorics, the related fundamental problem is to enumerate lattice animals of a specific size $n$; see \cite{bar3, bar4, BBLP,gut9, Viennot, conway1, mira} on this problem and references therein. Even for its simplest instance of enumeration of general animals, not much is known. Therefore, various constraints such as directional convexity and/or directional growth have been deployed to reduce the difficulty of the enumeration process in squared and hexagonal lattices.  See \cite{BBLP,buchin1, EB, ben1, melou1, DD, DV, kim1} for a few examples of the squared lattice and \cite{BBLP,Fer2, Fer1, GBL, Voge, lin1} for the hexagonal lattice, where several important subsets of animals including {\it bargraphs}, {\it column-convex animals,} and {\it convex animals} are considered. Although, a natural extension to these works is to explore the animals residing on the triangular lattice (see \cite{BBLP}), we were not able to find much work in this direction. Therefore, in this paper, we will strive to close some of the existing gaps concerning these objects.

{\it Polyiamonds}, which are animals residing on the triangular lattice, have been the subject of a few studies in recent years. For instance, in one of the first rigorous work on polyiamonds, Yang and Wilson \cite{yang1} proved the minimum perimeter of a polyiamond with $n$ cells is either $\ceil{\sqrt{6n}}$ or $\ceil{\sqrt{6n}}+1$ depending on the parity of $n$. To prove this result, they obtained a lower bound on the perimeter by considering maximal polyiamonds and showed how to construct minimal polyiamonds that attain the perimeter lower bounds. Extending in this direction, Malen and Rold\'an \cite{malen1} studied the maximum number of holes that a polyiamond with $n$ cells can enclose. Also, they obtained the minimum number of cells required to construct a polyiamond with $h$ holes. Very recently, in a series of works, by careful analyses of specific composition and concatenation operators on polyiamonds, Shahleh \cite{mira} and Barequet et al. \cite{bar1} improved the lower and upper bounds on the asymptotic growth rate with respect to the size to $2.8424$ and $3.6050$, respectively. Here, the growth rate of a sequence $a_n$ with respect to $n$ is simply $\lim_n a_n^{1/n}$ provided that this limit exists.

Our main results in this paper are related to the enumeration of several classes of polyiamonds; including baryiamonds, column-convex polyiamonds, and convex polyiamonds with respect to their perimeters. See Table \ref{tab-1} for an elaboration. This is mainly done by careful recursive analyses of these objects through column-by-column decomposition and cell pruning/growing arguments which lead to explicit forms of defined generating functions. We remark that one major difference between the triangular lattice and the squared and hexagonal lattices is the lack of some of the symmetrical structures in the former, which introduces an increasing level of complexity to the analysis. To that goal, we begin by providing a few definitions and notations.  Equip the two-dimensional plane with an $x-y$ system, where the angle between $x$-axis and $y$-axis is $\frac{\pi}{3}$ counter clock-wise. Mark each of these axes with integer points. We then partition this plane with equilateral triangles; each with a horizontal edge parallel to $x$-axis and the sides of length one.  We refer to these triangular building blocks as triangular cells or simply {\it cells}. The partition is uniquely defined such that there is one cell in the first quarter with a vertex at the intersection of the axes (i.e. origin) and two edges; coinciding with $x-y$ axes. We refer to this lattice, by $\calt.$ With some abuse of conventional notation, let the line $x=j$ denote the line parallel with the $y$-axis, intersection the $x$-axis at $j$. Define $y=j$ in a similar fashion. Considering the orientation of cells we identify two types of cells: (1) down-cells where their horizontal edges are at their bases, and (2) up-cells whose horizontal sides are on the top.

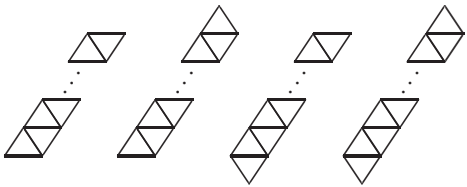
\begin{figure}[htp]
	\begin{picture}(55,30)
	\setlength{\unitlength}{.5mm}
	\linethickness{0.3mm}
	\def\uptt{\put(0,0){\line(1,0){10}}\put(10,0){\line(-2,3){5}}\put(0,0){\line(2,3){5}}}
	\def\dowtt{\put(0,10){\line(1,0){10}}\put(10,10){\line(-2,-3){5}}\put(0,10){\line(2,-3){5}}}
	\put(0,10){
		\put(0,0){\uptt\put(5,-2.5){\dowtt}}
		\put(5,7.5){\uptt\put(5,-2.5){\dowtt}}
		\multiput(16,17)(1.8,2.5){3}{\circle*{1}}
		\put(17,25){\uptt\put(5,-2.5){\dowtt}}
		\put(30,0){\put(0,0){\uptt\put(5,-2.5){\dowtt}}
			\put(5,7.5){\uptt\put(5,-2.5){\dowtt}}
			\multiput(16,17)(1.8,2.5){3}{\circle*{1}}
			\put(17,25){\uptt\put(5,-2.5){\dowtt}}
			\put(22,32.5){\uptt}}
		\put(60,0){\put(0,0){\uptt\put(5,-2.5){\dowtt}}
			\put(5,7.5){\uptt\put(5,-2.5){\dowtt}}
			\multiput(16,17)(1.8,2.5){3}{\circle*{1}}
			\put(17,25){\uptt\put(5,-2.5){\dowtt}}
			\put(0,-10){\dowtt}}
		\put(90,0){\put(0,0){\uptt\put(5,-2.5){\dowtt}}
			\put(5,7.5){\uptt\put(5,-2.5){\dowtt}}
			\multiput(16,17)(1.8,2.5){3}{\circle*{1}}
			\put(17,25){\uptt\put(5,-2.5){\dowtt}}
			\put(0,-10){\dowtt}\put(22,32.5){\uptt}}
	}
	\end{picture}
	\caption{From left to right, convex columns of type one, two, three, and four}\label{column-fig}
\end{figure}

\begin{figure}[htp]
	\begin{picture}(110,32)
	\setlength{\unitlength}{.5mm}
\def\uptt{\put(0,0){\line(1,0){10}}\put(10,0){\line(-2,3){5}}\put(0,0){\line(2,3){5}}} \def\dowtt{\put(0,10){\line(1,0){10}}\put(10,10){\line(-2,-3){5}}\put(0,10){\line(2,-3){5}}}
\put(0,15){
    \put(0,10){\put(0,0){\uptt\put(5,-2.5){\dowtt}}
			\put(5,7.5){\uptt\put(5,-2.5){\dowtt}}
			\put(10,15){\uptt\put(5,-2.5){\dowtt}}
			\put(0,-10){\dowtt}}
    \put(5,2.5){
        \put(0,0){\uptt\put(5,-2.5){\dowtt}}
		\put(5,7.5){\uptt\put(5,-2.5){\dowtt}}
		\put(10,15){\uptt\put(5,-2.5){\dowtt}}
        \put(15,22.5){\uptt\put(5,-2.5){\dowtt}}
        \put(20,30){\uptt\put(5,-2.5){\dowtt}}}
    \put(15,2.5){
        \put(0,0){\uptt\put(5,-2.5){\dowtt}}
		\put(5,7.5){\uptt\put(5,-2.5){\dowtt}}
		\put(10,15){\uptt\put(5,-2.5){\dowtt}}}
    \put(25,2.5){\put(0,0){\uptt\put(5,-2.5){\dowtt}}
			\put(5,7.5){\uptt}}
    \put(35,2.5){\put(0,0){\uptt\put(5,-2.5){\dowtt}}
            \put(5,7.5){\uptt\put(5,-2.5){\dowtt}}
			\put(10,15){\uptt}}
\multiput(-10,2.5)(0,-0.1){6}{\multiput(0,0)(8,0){8}{\line(1,0){5}}}
\multiput(-15,-12.5)(0.1,0){6}{\multiput(0,0)(8,12){5}{\line(2,3){5}}}
}
\put(80,15){
	\put(0,10){\put(0,0){\uptt\put(5,-2.5){\dowtt}}
			\put(5,7.5){\uptt\put(5,-2.5){\dowtt}}
			\put(10,15){\uptt\put(5,-2.5){\dowtt}}
			\put(0,-10){\dowtt}}
    \put(-5,-12.5){
        \put(0,0){\uptt\put(5,-2.5){\dowtt}}
		\put(5,7.5){\uptt\put(5,-2.5){\dowtt}}
		\put(10,15){\uptt\put(5,-2.5){\dowtt}}
        \put(15,22.5){\uptt\put(5,-2.5){\dowtt}}
        \put(20,30){\uptt\put(5,-2.5){\dowtt}}}
    \put(20,10){
        \put(0,0){\uptt\put(5,-2.5){\dowtt}}
		\put(5,7.5){\uptt\put(5,-2.5){\dowtt}}
		\put(10,15){\uptt\put(5,-2.5){\dowtt}}}
    \put(25,2.5){\put(0,0){\uptt\put(5,-2.5){\dowtt}}
			\put(5,7.5){\uptt}}
    \put(25,-12.5){\put(0,0){\uptt\put(5,-2.5){\dowtt}}
            \put(5,7.5){\uptt\put(5,-2.5){\dowtt}}
			\put(10,15){\uptt}}
\multiput(-10,2.5)(0,-0.1){6}{\multiput(0,0)(8,0){8}{\line(1,0){5}}}
\multiput(-15,-12.5)(0.1,0){6}{\multiput(0,0)(8,12){5}{\line(2,3){5}}}
}
\put(160,15){
	\put(0,10){\put(0,0){\uptt\put(5,-2.5){\dowtt}}
			\put(5,7.5){\uptt\put(5,-2.5){\dowtt}}
			\put(10,15){\uptt\put(5,-2.5){\dowtt}}
			\put(0,-10){\dowtt}}
    \put(-5,-12.5){
        \put(0,0){\uptt\put(5,-2.5){\dowtt}}
		\put(5,7.5){\uptt\put(5,-2.5){\dowtt}}
		\put(10,15){\uptt\put(5,-2.5){\dowtt}}
        \put(15,22.5){\uptt\put(5,-2.5){\dowtt}}
        \put(20,30){\uptt\put(5,-2.5){\dowtt}}}
    \put(5,-12.5){
        \put(0,0){\uptt\put(5,-2.5){\dowtt}}
			\put(5,7.5){\uptt\put(5,-2.5){\dowtt}}
			\put(10,15){\uptt\put(5,-2.5){\dowtt}}
			\put(0,-10){\dowtt}\put(15,22.5){\uptt}
			\put(5,7.5){\uptt}}
    \put(15,-12.5){\put(0,0){\uptt\put(5,-2.5){\dowtt}}
            \put(5,7.5){\uptt\put(5,-2.5){\dowtt}}
			\put(10,15){\uptt}}
    \put(25,-12.5){\put(0,0){\uptt\put(5,-2.5){\dowtt}}
            \put(5,7.5){\uptt\put(5,-2.5){\dowtt}}}
    \put(35,-12.5){\put(0,0){\uptt\put(5,-2.5){\dowtt}}
            \put(5,7.5){\uptt}}
    \put(45,-12.5){\put(0,0){\uptt\put(5,-2.5){\dowtt}}}
    \put(55,-12.5){\put(0,0){\uptt}}
\multiput(-10,2.5)(0,-0.1){6}{\multiput(0,0)(8,0){10}{\line(1,0){5}}}
\multiput(-15,-12.5)(0.1,0){6}{\multiput(0,0)(8,12){5}{\line(2,3){5}}}
}
\end{picture}
\caption{An example of a (left) baryiamond; (middle) column-convex polyiamond; (right) convex polyiamond}\label{fig2}
\end{figure}
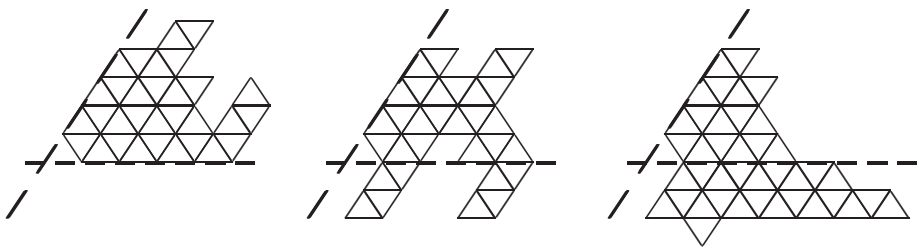

A {\it polyiamond} is a finite edge-connected set $\nu$ of cells in $\calt$, where for each distinct pair of cells $\Delta_1$ and $\Delta_2$ in $\nu$, there is a finite consecutive sequence of edge-connected cells in $\nu$ connecting $\Delta_1$ and $\Delta_2$. A cell in $\nu$ with at least one edge in common with a cell in $\nu^c:=\calt \setminus \nu$ is referred to as a boundary cell of $\nu$ and the common edge is referred to as a boundary edge. The {\it perimeter} of a polyiamond $\nu$ is  the total number of its boundary edges.  A {\it  column} (resp. {\it row}) of a polyiamond $\nu$ is a non-empty subset of cells in $\nu$, between lines $x=j$ and $x=j+1$ (resp. $y=j$ and $y=j+1$) for some $j\in\mathbb{Z}$. A polyiamond is called {\it column-convex} (resp. {\it row-convex}) if each of its columns (resp. rows) is a single contiguous block of cells. A {\it convex} polyiamond is both column-convex and row-convex. Depending on the orientation of the top and bottom cells of a convex column, we identify four types of convex columns as seen in Figure \ref{column-fig}.
A {\it baryiamond} is a column-convex polyiamond when either (a) it has only one convex column or (b) it has more than one column and the bottom cells of all consecutive columns are edge-connected. We use $\B,$ $\CCP,$, and $\CP$ to denote the sets of all unique baryiamonds, column-convex polyiamonds, and convex polyiamonds up-to translation in the triangular lattice. Clearly, $\B \subset \CCP$, $\CP \subset \CCP,$ and $\B \setminus \CP$ and $\CP \setminus \B$ are both non-empty.  See Figure \ref{fig2} for an example of elements in these sets.

The rest of the paper is organized as follows. Section \ref{bar-section} is dedicated to the analysis of $\B.$ In order to enumerate the number of baryiamonds of perimeter $n$, we study the generating function
\beqn \label{Bp-def}
B(p):=\sum_{\nu\in \B}p^{per(\nu)},
\feqn
where $per(\nu)$ denotes the perimeter of the polyiamond $\nu$. We obtain an explicit form for \eqref{Bp-def}. In particular, we prove
\begin{theorem} \label{bargraph-thm}
$$B(p)= \frac{1+p-2p^2-4p^3-3p^4+p^5-(1+p)^2\sqrt{(1-p)(1-p-2p^2-2p^3+p^4-p^5)}}{2p^2}.$$	Moreover, the number of baryiamonds with perimeter $n$ is asymptotic to
	$$\frac{(\xi+1)^2\sqrt{\xi^4+\xi^3-2\xi+1}}{2\sqrt{\pi n^3}}\xi^{-n-2},$$
	where $\xi=0.44617150675\cdots$ is the smallest root of the polynomial $1-p-2p^2-2p^3+p^4-p^5$.
\end{theorem}	

\par

In Section~\ref{CCP-section}, we focus on the column-convex polyiamonds, where we study the generating function
\beqn \label{Cp1-def}
C(p):=\sum_{\nu\in \CCP}p^{per(\nu)}.
\feqn
The main result of this section is

\begin{theorem} \label{ccp-thm}
	\beq
	&C(p^2)& = p^6+\frac{p^4(1-u_+^2)(1-u_-^2)A}{B}\\
&&=2p^6+3p^8+6p^{10}+15p^{12}+40p^{14}+113p^{16}+330p^{18}+988p^{20}+O(p^{22}),
    \feq
	where $u_+$ and $u_-$ are given in \eqref{equu1}-\eqref{equu2} and $A$ and $B$ are given in Lemma \ref{thCCPC1}.
	Moreover, the number of column-convex polyiamonds with perimeter $n$ is asymptotic to
	\beq
	\frac{(17997809\sqrt{17}+3^3\cdot13\cdot175463)\sqrt{95\sqrt{17}-119}}{2^7\cdot43^2\cdot 89^2\sqrt{6\pi n^3}}\left(\frac{3+\sqrt{17}}{2}\right)^{n-1}.
	\feq
\end{theorem}

Section~\ref{CP-section} is finally devoted to the study of convex polyiamonds and the generating function
\beqn \label{Cp-def}
F(p):=\sum_{\nu\in \CP}p^{per(\nu)}.
\feqn
We find the exact expression for this generating function and prove that

\begin{theorem}\label{mthcp}
	\begin{align*}
	&F(p)=-\frac{(2p+1)(p-1)^2(p+1)^3p^3}{(4p^2+2p-1)(4p^5-7p^3+p^2+4p-1)}\sqrt{1-2p-3p^2}\\
	&+\frac{(2p+1)(p-1)(p+1)^2p^5}{2(2p^4-p^3-2p^2+p+1)(5p^3+4p^2-p-1)(2p-1)}\sqrt{1-4p^2}\\
	&+\frac{1}{8}p^3-\frac{67}{40}p^2-\frac{103}{50}p-\frac{15071}{8000}+\frac{30+35p}{352(4p^2+2p-1)}\\
	&+\frac{318p^9+2083p^8+275p^7-3632p^6-1050p^5+2721p^4+680p^3-1093p^2-164p+169}{88(4p^{10}-8p^9+p^8+19p^7-8p^6-20p^5+9p^4+10p^3-5p^2-2p+1)}\\
	&+\frac{-238p^4+3p^3+537p^2+214p-149}{704(4p^5-7p^3+p^2+4p-1)}\\
	&+\frac{119p^2+114p+19}{500(5p^3+4p^2-p-1)}+
	\frac{3p^3+2p^2-p-1}{8(2p^4-p^3-2p^2+p+1)}\\
	&=2p^3+3p^4+6p^5+15p^6+38p^7+102p^8+272p^9+739p^{10}+2006p^{11}+O(p^{12}).
	\end{align*}
	In addition, the number of convex polyiamonds with perimeter $n$ is asymptotic to
	$$\frac{1280}{441\sqrt{3\pi n^3}}3^n.$$
\end{theorem}

At last, we summarized our findings along with known results for other lattices in Table \ref{tab-1} for the sake of comparison. We remark that most of the calculations in this paper are completed by extensive use of {\it Maple} and we chose to not include some of the details for the sake of brevity and space.

\begin{table}
\centering \tiny
\begin{tabular}{|l|l|l|l|}\hline
Lattice& Animals sub-family& The number of animals in sub-family with perimeter $n$ &Reference\\[6pt]\hline\hline
Triangular	& Bargraphs (Baryiamonds) &$\frac{(\xi+1)^2\sqrt{\xi^4+\xi^3-2\xi+1}}{2\sqrt{\pi n^3}}\xi^{-n-2}$&Theorem \ref{bargraph-thm}\\[4pt]\cline{2-4}
& & &\\[-6pt]
& Column-Convex Animals&	$\frac{(17997809\sqrt{17}+3^3\cdot13\cdot175463)\sqrt{95\sqrt{17}-119}}{2^7\cdot43^2\cdot 89^2\sqrt{6\pi n^3}}\left(\frac{3+\sqrt{17}}{2}\right)^{n-1}$& Theorem \ref{ccp-thm}\\[4pt]\cline{2-4}
& & &\\[-6pt]
&Convex Animals& $\frac{1280}{441\sqrt{3\pi n^3}}3^n$ & Theorem \ref{mthcp} \\[4pt]
\hline\hline
& & &\\[-6pt]
Squared	& Bargraphs&$\sqrt{\frac{2(1-\rho-\rho^3)}{\pi\rho n^3}}\rho^{-n/2}$&\\
&&with $\rho=\frac{1}{3}\left(-1-\frac{2^{8/3}}{(13+3\sqrt{3})^{1/3}}+2^{1/3}(13+3\sqrt{33})^{1/3}\right)$&\cite{Pen1}\\\cline{2-4}
& & &\\[-6pt]
 &Column-Convex Animals& $\frac{2\sqrt{5\sqrt{2}-7}((138\sqrt{2}+444)\sqrt{10-7\sqrt{2}}
+589\sqrt{2}-410)}{2209}\cdot
\frac{(1+\sqrt{2})^{n}}{n\sqrt{\pi n/2}}$&\cite{Fer1}\\[6pt]\cline{2-4}
& & &\\[-6pt]
&Convex Animals&  $(n+3)2^{n-8}-4(n-7){n-8 \choose n/2-4}$&\cite{DV}\\[4pt]
\hline\hline
Hexagonal & Bargraphs & $\frac{\sqrt{3}(13-\sqrt{17})}{\sqrt{17\pi n^3}}\sqrt{2}^{\,n}$ & \cite{toufik10}\\[6pt]\cline{2-4}
& Column-Convex Animals&$\frac{32(179\sqrt{7}-140)\sqrt{2}}{22743\sqrt{3\pi n^3}}\sqrt{3}^n$& g.f. reported in \cite{Fer1}\\[6pt]\cline{2-4}
& Convex Animals& Unknown & g.f. reported in \cite{GBL}\\[6pt]\hline
\end{tabular}
\caption{\label{tab-1} Summary of enumeration of bargraphs, column-convex animals, and convex animals on different lattices with respect to the perimeter. The reported results for convex animals on the squared lattice are exact numbers, whereas all the others are asymptotic results.}
\end{table}

\section{Baryiamonds} \label{bar-section}
\begin{figure}[htp]
	\begin{picture}(95,27)
	\setlength{\unitlength}{.5mm}
	\def\uptt{\put(0,0){\line(1,0){10}}\put(10,0){\line(-2,3){5}}\put(0,0){\line(2,3){5}}} \def\dowtt{\put(0,10){\line(1,0){10}}\put(10,10){\line(-2,-3){5}}\put(0,10){\line(2,-3){5}}}
	\put(0,0){
		\put(0,0){\uptt\put(5,-2.5){\dowtt}}
		\put(5,7.5){\uptt\put(5,-2.5){\dowtt}}
		\multiput(16,17)(1.8,2.5){3}{\circle*{1}}
		\put(17,25){\uptt\put(5,-2.5){\dowtt}}
		\put(22,32.5){\uptt\put(5,-2.5){\dowtt}}
		\put(25,0){
			\put(0,0){\uptt\put(5,-2.5){\dowtt}}
			\put(5,7.5){\uptt\put(5,-2.5){\dowtt}}
			\multiput(16,17)(1.8,2.5){3}{\circle*{1}}
			\put(17,25){\uptt\put(5,-2.5){\dowtt}}
			\put(22,32.5){\uptt\put(5,-2.5){\dowtt}}}
		\put(35,0){
			\put(0,0){\uptt\put(5,-2.5){\dowtt}}
			\multiput(11,10.5)(1.8,2.5){3}{\circle*{1}}
			\put(12,18){\uptt\put(5,-2.5){\dowtt}}}
		\put(65,0){
			\put(0,0){\uptt\put(5,-2.5){\dowtt}}
			\put(5,7.5){\uptt\put(5,-2.5){\dowtt}}
			\multiput(16,17)(1.8,2.5){3}{\circle*{1}}
			\put(17,25){\uptt\put(5,-2.5){\dowtt}}
			\put(22,32.5){\uptt\put(5,-2.5){\dowtt}}}
		\put(75,0){
			\put(0,0){\uptt\put(5,-2.5){\dowtt}}
			\multiput(11,10.5)(1.8,2.5){3}{\circle*{1}}
			\put(12,18){\uptt\put(5,-2.5){\dowtt}}
			\put(17,25.5){\uptt}}
		\put(125,0){
			\put(0,0){\uptt\put(5,-2.5){\dowtt}}
			\put(5,7.5){\uptt\put(5,-2.5){\dowtt}}
			\multiput(16,17)(1.8,2.5){3}{\circle*{1}}
			\put(17,25){\uptt\put(5,-2.5){\dowtt}}
			\put(22,32.5){\uptt\put(5,-2.5){\dowtt}}}
		\put(135,0){
			\put(0,0){\uptt\put(5,-2.5){\dowtt}}
			\multiput(11,10.5)(1.8,2.5){3}{\circle*{1}}
			\multiput(37,49)(1.8,2.5){3}{\circle*{1}}
			\put(12,17.5){\uptt\put(5,-2.5){\dowtt}}
			\put(17,25){\uptt\put(5,-2.5){\dowtt}}
			\put(22,32.5){\uptt\put(5,-2.5){\dowtt}}
			\put(27,40){\uptt\put(5,-2.5){\dowtt}}}
		\put(165,0){
			\put(0,0){\uptt\put(5,-2.5){\dowtt}}
			\put(5,7.5){\uptt\put(5,-2.5){\dowtt}}
			\multiput(16,17)(1.8,2.5){3}{\circle*{1}}
			\multiput(46,47)(1.8,2.5){3}{\circle*{1}}
			\put(17,25){\uptt\put(5,-2.5){\dowtt}}
			\put(22,32.5){\uptt\put(5,-2.5){\dowtt}}}
		\put(175,0){
			\put(0,0){\uptt\put(5,-2.5){\dowtt}}
			\multiput(11,10.5)(1.8,2.5){3}{\circle*{1}}
			\put(12,17.5){\uptt\put(5,-2.5){\dowtt}}
			\put(17,25){\uptt\put(5,-2.5){\dowtt}}
			\put(22,32.5){\uptt\put(5,-2.5){\dowtt}}
			\put(27,40){\uptt}}
	}
	\end{picture}
	\caption{Various possibilities of baryiamonds in $\B_k^{(1)}$, according to the type of their second column (if any)}\label{figBB1}
\end{figure}
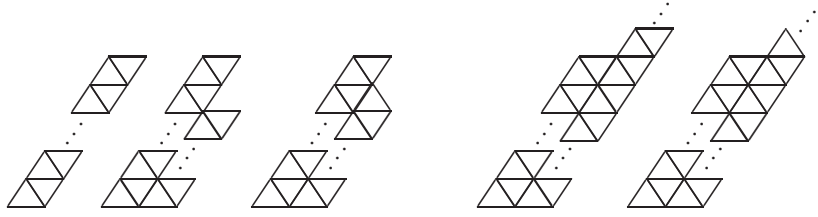

In this section, we enumerate the baryiamonds with respect to their perimeter. Recall the types of convex columns in Figure \ref{column-fig}. Let $\B_k^{(i)}$ be the set of all baryiamonds in $\B$ whose first columns are of type $i$ and have $k$ up-cells. Set $\B^{(i)}=\cup_{k\geq 1} \B_k^{(i)}.$ We define $B^{(i)}_k(p)$ to be the generating function with respect to perimeter over $\B_k^{(i)}$; that is,
\beq
B^{(i)}_k(p)= \sum_{\nu\in\B_k^{(i)}} p^{per(\nu)},
\feq
where $per(\nu)$ is the perimeter of $\nu.$ We similarly define the generating function with respect to perimeter over $\B^{(i)}$ for $i=1,2,3,4$:
\beq
B^{(i)}(p,u)= \sum_{k\geq 1} \sum_{\nu\in\B_k^{(i)}} u^k p^{per(\nu)},
\feq
where $u$ marks the number of up-cells in the first column.   It is clear that $B(p)$, as defined in \eqref{Bp-def}, is $\sum_{i=1}^4 B^{(i)}(p,1).$ In the next few paragraphs, we will obtain explicit expressions for $B^{(i)}(p,u)$ and prove Theorem \eqref{bargraph-thm}. The type of methods used to derive these expressions are referred to {\it cell growing/pruning} methods throughout this paper. First, note that each $\nu\in \B_k^{(1)}$ falls exclusively into one of the cases described in Figure \ref{figBB1}. Adding up the contribution of each of these cases, we arrive at
\begin{align} \label{bar1eq}
B_k^{(1)}(p)=p^{2k+2}&+\sum_{j=1}^kp^{2k+2-2j}(B^{(1)}_j(p)+B^{(2)}_{j-1}(p))\notag\\
&+\sum_{j\geq k+1}p^2(B^{(1)}_j(p)+B^{(2)}_{j-1}(p)).
\end{align}

Next, we describe $B^{(i)}_k(p)$, $i=2,3,4$, with respect to $B^{(1)}_k(p)$. To that goal, note that
\begin{itemize}
\item[(i)] By removing the top down-cell of the first column in $\nu\in \B_k^{(2)}$ with more than one cell, we get a bargraph in $\B_k^{(1)}$; and hence for all $k\geq1$,
\beqn \label{bar2eq}
B^{(2)}_k(p)=pB^{(1)}_k(p);
\feqn
\item[(ii)] By adding a down-cell to the bottom of the first column of $\B_k^{(3)}$, we get a bargraph in $\B_k^{(1)}$; hence, for all $k\geq1$,
\beqn \label{bar3eq}
B^{(3)}_k(p)=\frac{1}{p}B^{(1)}_k(p);
\feqn
\item[(iii)] By removing the top down-cell of the first column and adding a down-cell to the bottom of the same column of the $\nu \in \B_k^{(4)}$, we see that, for all $k\geq1,$
\beqn \label{bar4eq}
B^{(4)}_k(p)=B^{(1)}_k(p).
\feqn
\item[(iv)] The initial conditions are given by
\beqn \label{bar5eq}
B_0^{(2)}(p)=p^3,  \quad B^{(i)}_0(p)=0 \quad i=1,3,4.
\feqn
\end{itemize}
By multiplying each of equations in \eqref{bar1eq}-\eqref{bar5eq} by $u^k$ and adding up over $k\geq1$, we get
\beqn \label{Biup-bar}
&&B^{(1)}(p,u)=\frac{p^4u}{1-p^2u}+\frac{p^2}{1-p^2u}(B^{(1)}(p,u)+uB^{(2)}(p,u))\notag \\
&&\qquad\qquad\quad  +\frac{p^2}{1-u}(uB^{(1)}(p,1)+uB^{(2)}(p,1)-B^{(1)}(p,u)-uB^{(2)}(p,u)), \notag\\
&&B^{(2)}(p,u)=p^3+pB^{(1)}(p,u), \quad B^{(3)}(p,u)=\frac{1}{p}B^{(1)}(p,u), \notag \\
&& B^{(4)}(p,u)=B^{(1)}(p,u).
\feqn
This implies
\begin{align}
B^{(1)}(p,u)&=\frac{p^4u}{1-p^2u}+\frac{p^2}{1-p^2u}((1+pu)B^{(1)}(p,u)+p^3u)\nonumber\\
&+\frac{p^2}{1-u}(u(1+p)B^{(1)}(p,1)-(1+pu)B^{(1)}(p,u)).\label{eqBB11}
\end{align}
This functional equation is then solved by the kernel method. See \cite{Ban} for an introduction to this method. To do so, we re-arrange the equation such that all the terms containing $B^{(1)}(p,u)$ are on the left side and the rest of terms are collected on the right side. Setting $u$ to  $$u_0=\frac{1+p^4-(1+p)\sqrt{(1-p)(1-p-2p^2-2p^3+p^4-p^5)}}{2p^2(1+p-p^3)}$$
will diminish the left hand side, and therefore, we obtain
\beqn \label{B1p1}
B^{(1)}(p,1)&=\frac{1-p-p^2-p^3-\sqrt{(1-p)(1-p-2p^2-2p^3+p^4-p^5)}}{2p}.
\feqn
Then, by substituting \eqref{B1p1} into the last three equations of \eqref{Biup-bar} and adding them up, we complete the proof of the first part of Theorem \ref{bargraph-thm}. The second part is easily obtained by considering the singularity analysis (for example, see Chapter 6 in \cite{FS}) applied to the generating function $B(p)$.

\section{Column-convex polyiamonds} \label{CCP-section}

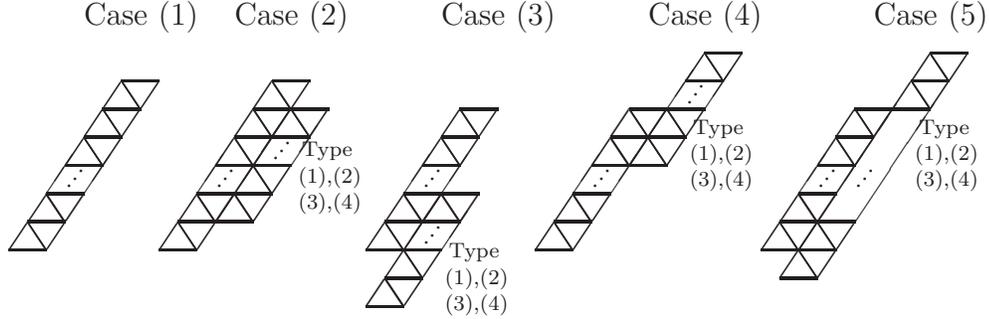
\begin{figure}[htp]
\begin{picture}(90,35)
\setlength{\unitlength}{.5mm}
\linethickness{0.3mm}
\put(-20,70){Case (1)}
\put(20,70){Case (2)}\put(75,70){Case (3)}\put(130,70){Case (4)}
\put(190,70){Case (5)}
\def\uptt{\put(0,0){\line(1,0){10}}\put(10,0){\line(-2,3){5}}\put(0,0){\line(2,3){5}}}
\def\dowtt{\put(0,10){\line(1,0){10}}\put(10,10){\line(-2,-3){5}}\put(0,10){\line(2,-3){5}}}
\put(-40,10){\put(0,0){\uptt\put(5,-2.5){\dowtt}}
\put(5,7.5){\uptt\put(5,-2.5){\dowtt}}
\multiput(16,17)(1.6,2){3}{\circle*{1}}
\put(15,22.5){\uptt\put(5,-2.5){\dowtt}}
\put(20,30){\uptt\put(5,-2.5){\dowtt}}
\put(25,37.5){\uptt\put(5,-2.5){\dowtt}}
\put(10,15){\line(2,3){7}}\put(20,15){\line(2,3){7}}}
\put(0,10){\put(0,0){\uptt\put(5,-2.5){\dowtt}}
\put(5,7.5){\uptt\put(5,-2.5){\dowtt}}
\multiput(16,17)(1.6,2){3}{\circle*{1}}
\put(15,22.5){\uptt\put(5,-2.5){\dowtt}}
\put(20,30){\uptt\put(5,-2.5){\dowtt}}
\put(25,37.5){\uptt\put(5,-2.5){\dowtt}}
\put(10,15){\line(2,3){7}}\put(20,15){\line(2,3){7}}}
\put(15,17.5){\put(0,0){\uptt\put(5,-2.5){\dowtt}}
\put(5,7.5){\uptt\put(5,-2.5){\dowtt}}
\multiput(16,17)(1.6,2){3}{\circle*{1}}
\put(15,22.5){\uptt\put(5,-2.5){\dowtt}}
\put(10,15){\line(2,3){7}}\put(20,15){\line(2,3){7}}}
\put(38,35){\tiny Type}
\put(37,28){\tiny (1),(2)}\put(37,21){\tiny (3),(4)}
\put(55,0){
\put(0,10){\put(0,0){\uptt\put(5,-2.5){\dowtt}}
\put(5,7.5){\uptt\put(5,-2.5){\dowtt}}
\multiput(16,17)(1.6,2){3}{\circle*{1}}
\put(15,22.5){\uptt\put(5,-2.5){\dowtt}}
\put(20,30){\uptt\put(5,-2.5){\dowtt}}
\put(10,15){\line(2,3){7}}\put(20,15){\line(2,3){7}}}
\put(-15,-22.5){\put(15,17.5){\put(0,0){\uptt\put(5,-2.5){\dowtt}}
\put(5,7.5){\uptt\put(5,-2.5){\dowtt}}
\multiput(16,17)(1.6,2){3}{\circle*{1}}
\put(15,22.5){\uptt\put(5,-2.5){\dowtt}}
\put(10,15){\line(2,3){7}}\put(20,15){\line(2,3){7}}}}
\put(22,8){\tiny Type}\put(21,1){\tiny (1),(2)}\put(21,-6){\tiny (3),(4)}
}
\put(100,0){
\put(0,10){\put(0,0){\uptt\put(5,-2.5){\dowtt}}
\put(5,7.5){\uptt\put(5,-2.5){\dowtt}}
\multiput(16,17)(1.6,2){3}{\circle*{1}}
\put(15,22.5){\uptt\put(5,-2.5){\dowtt}}
\put(20,30){\uptt\put(5,-2.5){\dowtt}}
\put(10,15){\line(2,3){7}}\put(20,15){\line(2,3){7}}}
\put(10,15){\put(15,17.5){\put(0,0){\uptt\put(5,-2.5){\dowtt}}
\put(5,7.5){\uptt\put(5,-2.5){\dowtt}}
\multiput(16,17)(1.6,2){3}{\circle*{1}}
\put(15,22.5){\uptt\put(5,-2.5){\dowtt}}
\put(10,15){\line(2,3){7}}\put(20,15){\line(2,3){7}}}}
\put(20,33){\put(22,8){\tiny Type}\put(21,1){\tiny (1),(2)}\put(21,-6){\tiny (3),(4)}}
}
\put(160,0){
\put(0,10){\put(0,0){\uptt\put(5,-2.5){\dowtt}}
\put(5,7.5){\uptt\put(5,-2.5){\dowtt}}
\multiput(16,17)(1.6,2){3}{\circle*{1}}
\put(15,22.5){\uptt\put(5,-2.5){\dowtt}}
\put(20,30){\uptt\put(5,-2.5){\dowtt}}
\put(10,15){\line(2,3){7}}\put(20,15){\line(2,3){7}}}
\put(-10,-15){\put(15,17.5){\put(0,0){\uptt\put(5,-2.5){\dowtt}}
\put(5,7.5){\uptt\put(5,-2.5){\dowtt}}
\multiput(21,24.5)(1.6,2){3}{\circle*{1}}
\put(30,45){\uptt\put(5,-2.5){\dowtt}}
\put(35,52.5){\uptt\put(5,-2.5){\dowtt}}
\put(10,15){\line(2,3){21}}\put(20,15){\line(2,3){21}}}}
\put(20,33){\put(22,8){\tiny Type}\put(21,1){\tiny (1),(2)}\put(21,-6){\tiny (3),(4)}}
}
\end{picture}
\caption{Various possibilities of polyiamonds in $\CCP_k^{(1)}$ with respect to the relative position of the first and second columns; Only the cases whose second columns are of type $1$ are shown (See Figure \ref{column-fig}). The second column can be replaced with any other types mentioned in the picture.}\label{ccp2}
\end{figure}

This section is devoted to the enumeration of column-convex polyiamonds with respect to their perimeter. Let $\CCP_k^{(i)}$ be the set of all polyiamonds in $\CCP$ where their first columns are of type $i$ and have $k$ up-cells.
We define $C^{(i)}_k(p)$ to be the generating function with respect to perimeter over $\CCP_k^{(i)}$; that is,
\beq
C^{(i)}_k(p)= \sum_{\nu \in \CCP_k^{(i)}} p^{per(\nu)}.
\feq
Moreover, for $i=1,2,3,4$, we set $$C^{(i)}(p,u)=\sum_{k\geq0}C_k^{(i)}(p)u^{k}.$$
With these notations, $C(p)$ as defined by \eqref{Cp1-def} is exactly $\sum_{i=1}^4 C^{(i)}(p,1).$ Therefore, we focus on obtaining $C^{(i)}(p,1)$ for $i=1,2,3,4$.
First, let us review a few observations that relate $C_k^{(i)}(p),$ $i=2,3,4,$ to $C_k^{(1)}(p).$
\begin{itemize}
\item[(i)] By removing the top down-cell of the first column of $\nu\in \CCP_k^{(2)}$, we get another column-convex polyiamond in $\CCP_k^{(1)}$; and hence for all $k\geq1,$ we have
\beqn \label{eq1CCP}
C^{(2)}_k(p)=pC^{(1)}_k(p);
\feqn
\item[(ii)] By adding a down-cell to the bottom of the first column of $\nu \in \CCP_k^{(3)}$, we get another column-convex polyiamond in $\CCP_k^{(1)}$; hence, for all $k\geq1$
\beqn \label{eq2CCP}
C^{(3)}_k(p)=\frac{1}{p}C^{(1)}_k(p);
\feqn
\item[(iii)] By removing the top down-cell of the first column and adding a down-cell to the bottom of the same column of $\nu \in \CCP_k^{(4)}$, we see that, for all $k\geq1,$
\beqn \label{eq3CCP}
C^{(4)}_k(p)=C^{(1)}_k(p);
\feqn
\item[(iv)] The initial conditions are given by
\beqn \label{eq4CCP}
C^{(2)}_0(p)=p^3, \quad C^{(j)}_0(p)=0, \quad j=1,3,4.
\feqn
\end{itemize}
Multiplying the equations \eqref{eq1CCP}-\eqref{eq4CCP} by $u^k$ and adding all the terms up over $k\geq 0$, we arrive at the system of equations
\beqn\label{eqC234}
C^{(2)}(p,u)=p^3 +pC^{(1)}(p,u),\quad C^{(3)}(p,u)=\frac{1}{p}C^{(1)}(p,u),\quad C^{(4)}(p,u)=C^{(1)}(p,u).
\feqn
Hence, to calculate $C(p),$ it is enough to find the explicit form of $C^{(1)}(p,u).$ To that goal, we a write a recurrence for $C^{(1)}_k(p)$ by investigating the types of first and second (if any) columns and their relative positions to each other.  See Figure~\ref{ccp2} for a few examples of possible cases. Let $\nu\in \CCP^{(1)}_k$ with $k\geq 2$.

Adding up the contributions of all cases in Figure \ref{ccp2}, we obtain
\begin{align*}
C_k^{(1)}(p)=p^{2k+2}
&+\sum_{j=1}^k(k+1-j)p^{2k+2-2j}(C_j^{(1)}(p)+C_{j-1}^{(2)}(p)+C_{j+1}^{(3)}(p)+C_{j}^{(4)}(p))\\
&+2\sum_{j=2}^k\frac{p^{2k+4-2j}-p^{2k+2}}{1-p^2}(C_j^{(1)}(p)+C_{j-1}^{(2)}(p)+C_{j+1}^{(3)}(p)+C_{j}^{(4)}(p))\\
&+2\sum_{j\geq k+1}\frac{p^2-p^{2k+2}}{1-p^2}(C_j^{(1)}(p)+C_{j-1}^{(2)}(p)+C_{j+1}^{(3)}(p)+C_{j}^{(4)}(p))\\
&+\sum_{j\geq k+2}(j-k-1)p^2(C_j^{(1)}(p)+C_{j-1}^{(2)}(p)+C_{j+1}^{(3)}(p)+C_{j}^{(4)}(p)),
\end{align*}
for all $k\geq2$.

By multiplying by $u^k$ and summing over $k\geq2$, we then obtain
\begin{align*}
C^{(1)}(p,u)&=\frac{p^4u}{1-p^2u}\\
&+\frac{p^2}{(1-p^2u)^2}(C^{(1)}(p,u)+uC^{(2)}(p,u)+\frac{1}{u}(C^{(3)}(p,u)-uC^{(3)}_1(p))+C^{(4)}(p,u))\\
&-\frac{2p^2}{(1-p^2u)(1-u)}(C^{(1)}(p,u)+uC^{(2)}(p,u)+\frac{1}{u}(C^{(3)}(p,u)-uC^{(3)}_1(p))+C^{(4)}(p,u))\\
&+\frac{2p^2u}{(1-u)(1-p^2u)}(C^{(1)}(p,1)+C^{(2)}(p,1)+C^{(3)}(p,1)-C^{(3)}_1(p)+C^{(4)}(p,1))\\
&+\frac{p^2}{(1-u)^2}(C^{(1)}(p,u)+uC^{(2)}(p,u)+\frac{1}{u}(C^{(3)}(p,u)-uC^{(3)}_1(p))+C^{(4)}(p,u))\\
&-\frac{p^2u(2-u)}{(1-u)^2}(C^{(1)}(p,1)+C^{(2)}(p,1)+C^{(3)}(p,1)-C^{(3)}_1(p)+C^{(4)}(p,1))\\
&+\frac{p^2u}{1-u}\frac{\partial}{\partial u}(C^{(1)}(p,u)+uC^{(2)}(p,u)+\frac{1}{u}(C^{(3)}(p,u)-uC^{(3)}_1(p))+C^{(4)}(p,u))\mid_{u=1}.
\end{align*}
To compute $C^{(1)}(p^2,1),$ we next define $D^{(i)}(p,u):=C^{(i)}(p^2,u^2)$. By \eqref{eqC234}, the last equation can be written as
\begin{align}\label{CCPeqL1}
&\left(1-\frac{p^2u^2(1+p^2u^2)^2(1-p^4)^2}{(1-p^4u^2)^2(1-u^2)^2}\right)D^{(1)}(p,u)\nonumber\\
&=\frac{p^2u^2(1+p^2)(p^6u^2-2p^4u^4+3p^4u^2-2p^2u^2+p^2-1)}{(1-p^4u^2)(1-u^2)^2}D^{(1)}(p,1)\\
&-\frac{p^{10}u^4}{(1-p^4u^2)^2}C^{(1)}_1(p^2)+\frac{p^2u^2(1+p^2)^2}{1-u^2}
\frac{\partial}{\partial u}C^{(1)}(p^2,u)\mid_{u=1}
-\frac{(p^4u^2-p^2-1)p^8u^2}{(1-p^4u^2)^2}.\nonumber
\end{align}
Recall that $D^{(1)}(p,0)=C^{(1)}(p^2,0)=0.$ We take the derivative at $u=0$, and solve for $C^{(1)}_1(p^2)$. This leads to
\beqn \label{C1p2ueq1}
C^{(1)}_1(p^2)=(1+p)p^4+p(p^2-1)D^{(1)}(p,1) +p(p+1)^2\frac{\partial}{\partial u}C^{(1)}(p^2,u)\mid_{u=1}.
\feqn
Then, substituting \eqref{C1p2ueq1} into \eqref{CCPeqL1}, we arrive at
\begin{align}\label{CCPeqL2}
&\left(1-\frac{p^2u^2(1+p^2u^2)^2(1-p^4)^2}{(1-p^4u^2)^2(1-u^2)^2}\right)D^{(1)}(p,u)\nonumber\\
&=\biggl(\frac{p^2u^2(1+p^2)(1+(2u^2-1)p^2+2u^2(u^2-2)p^4-2u^4p^6)}
{(1-p^4u^2)^2(1-u^2)^2}\nonumber\\
&\qquad+\frac{p^{10}u^4(1+p^2)(u^2(3-2u^2)-(u^4-3u^2+1)p^{2}+(1-u^2)^2p^{4})}
{(1-p^4u^2)^2(1-u^2)^2}\biggr)D^{(1)}(p,1)\\
&+\frac{p^2u^2(1+p^2)^2(1-2u^2p^4+u^4p^8-u^2(1-u^2)p^{10})}{(1-p^4u^2)^2(1-u^2)}\frac{\partial}{\partial u}C^{(1)}(p^2,u)\mid_{u=1}\nonumber\\
&+\frac{p^8u^2(1+p^2-u^2p^4-u^2p^{10}-u^2p^{12})}{(1-p^4u^2)^2}.\nonumber
\end{align}
Again, we deploy the kernel method to solve this functional equation. To that end, let $K(u)=1-\frac{p^2u^2(1+p^2u^2)^2(1-p^4)^2}{(1-p^4u^2)^2(1-u^2)^2}$ and set
\begin{align}
&u_+^2=u_+^2(p)=\frac{p^{10}-2p^6+2p^4+p^2+2}{4p^4}
-\frac{(1-p^4)(1+p^2)}{4p^3}\sqrt{1+p^2}\sqrt{4-3p^2+p^4}\label{equu1}\\
&-\frac{1-p^4}{2p^4}
\sqrt{1+p^2(1+p^2)^2+\frac{p(p^{10}-2p^6+2p^4+p^2+2)}{2(1-p^2)}(p-p^3-\sqrt{1+p^2}\sqrt{4-3p^2+p^4})}\nonumber\\
&=1+p+\frac{1}{2}p^2+\frac{9}{8}p^3+2p^4+\frac{239}{128}p^5+O(p^6),\nonumber\\
&u_-^2=u_-^2(p)=u_+^2(-p),\label{equu2}
\end{align}
where $u_{\pm}$ are the roots of $K(u)=0$. Therefore,
 \begin{align}\label{CCPeqL3}
&\biggl(1+(2u_\pm^2-1)p^2+2u_\pm^2(u_\pm^2-2)p^4-2u_\pm^4p^6\nonumber\\
&\qquad+p^{8}u_\pm^2(u_\pm^2(3-2u_\pm^2)-(u_\pm^4-3u_\pm^2+1)p^{2}+(1-u_\pm^2)^2p^{4})\biggr)D^{(1)}(p,1)\nonumber\\
&=(1+p^2)(1-u_\pm^2)(1-2u_\pm^2p^4+u_\pm^4p^8-u_\pm^2(1-u_\pm^2)p^{10})\frac{\partial}{\partial u_\pm}C^{(1)}(p^2,u_\pm)\mid_{u_\pm=1}\nonumber\\
&+\frac{p^6(1-u_\pm^2)^2(1+p^2-u_\pm^2p^4-u_\pm^2p^{10}-u_\pm^2p^{12})}{1+p^2}.\nonumber
\end{align}
Solving this system of equations, we obtain an explicit formula for $C^{(1)}(p^2,1)$. More specifically,
\begin{lemma}\label{thCCPC1}
The generating function $C^{(1)}(p^2,1)$ is given by
\begin{align*}
\frac{p^6(1-u_+^2)(1-u_-^2)A}{(1+p^2)^2B},
\end{align*}
where
\begin{align*}
A&=2p^6+p^4-p^2-1-p^4(p^6+p^4-1)(u_+^2+u_-^2)\\
&\qquad\qquad+p^8(p^{12}+p^{10}-2p^8-2p^6+2p^4+2p^2-1)u_+^2u_-^2,\\
B&=1-2p^4+2p^8(u_+^2+u_-^2)-p^8(u_+^4+u_-^4)
-2p^4(p^{14}-2p^{10}+p^8+p^6+1)u_+^2u_-^2\\
&\qquad\qquad+p^8(p^{10}-2p^6+p^2+2)u_+^2u_-^2(u_+^2+u_-^2)+p^{12}(p^4-2)u_+^4u_-^4.
\end{align*}
\end{lemma}
This lemma along with \eqref{eqC234} completes the proof of the first part in Theorem \ref{ccp-thm}. To show the second part of Theorem \ref{ccp-thm}, we apply the singularity analysis to the generating function $C(p^2)$. Note that $C(p^2)$ can be further rewritten as
{\small$$\frac{a_{00}+a_{10}\sqrt{\alpha}+a_{01}\sqrt{\beta}+a_{11}\sqrt{\alpha\beta}}
{2(4p^{22}-16p^{20}+20p^{18}+32p^{16}-75p^{14}-76p^{12}+182p^{10}+152p^8-155p^6-164p^4+48p^2+64)},$$}
where
\begin{align*}
a_{00}&=(1-p^2)\bigl(14p^{24}-34p^{22}+3p^{20}+143p^{18}-21p^{16}-377p^{14}-57p^{12}+449p^{10}\\
&+181p^8-245p^6-112p^4+72p^2+32\bigr),\\
a_{10}&=\frac{\sqrt{2}}{4}(p^4-1)\bigl(4p^{24}-12p^{22}+2p^{20}+54p^{18}-19p^{16}-138p^{14}+54p^{12}+250p^{10}\\
&-3p^8-234p^6-62p^4+80p^2+32\bigr)\\
&+\frac{\sqrt{2}p}{4}(p^2+1)(p^2-1)^2\bigl(4p^{18}-4p^{16}-10p^{14}+30p^{12}+63p^{10}-22p^8-99p^6\\
&-28p^4+56p^2+32\bigr)\sqrt{p^2(p^2-1)^2+4},
\end{align*}
\begin{align*}
a_{01}&=\frac{\sqrt{2}}{4}(p^4-1)\bigl(4p^{24}-12p^{22}+2p^{20}+54p^{18}-19p^{16}-138p^{14}+54p^{12}+250p^{10}\\
&-3p^8-234p^6-62p^4+80p^2+32\bigr)\\
&-\frac{\sqrt{2}p}{4}(p^2+1)(p^2-1)^2\bigl(4p^{18}-4p^{16}-10p^{14}+30p^{12}+63p^{10}-22p^8-99p^6\\
&-28p^4+56p^2+32\bigr)\sqrt{p^2(p^2-1)^2+4},\\
a_{11}&=p^2(p^2+1)(p^2-1)^2(2p^{12}+p^{10}+4p^6+5p^4-3p^2-4), \\
\alpha&=p^{12}-2p^8+4p^6+5p^4+4p^2+2+\frac{p(p^{10}-2p^6+2p^4+p^2+2)}{1-p^2}\sqrt{p^2(p^2-1)^2+4},\\
\beta&=p^{12}-2p^8+4p^6+5p^4+4p^2+2-\frac{p(p^{10}-2p^6+2p^4+p^2+2)}{1-p^2}\sqrt{p^2(p^2-1)^2+4}.
\end{align*}
Consider $r:=\frac{\sqrt{\sqrt{17}-3}}{2}$. Then, $\alpha$ and $\beta$ can be written as
\begin{align*}
\alpha^2&=\frac{95\sqrt{17}-119}{6}r^4(1+p/r)+O((1+p/r)^2),\\
\beta^2&=\frac{95\sqrt{17}-119}{6}r^4(1-p/r)+O((1-p/r)^2).
\end{align*}
Consequently, when $p^2$ is near the dominant singularity of $C(p^2)$,  $C(p^2)$ can be approximated by
\beq
 -\frac{(17997809\sqrt{17}+61587513)\sqrt{2}}{937339456}\frac{\sqrt{95\sqrt{17}-119}r^2}{\sqrt{6}}
\left(\sqrt{1+p/r}+\sqrt{1-p/r}\right),
\feq
from which the second part of Theorem \ref{ccp-thm} follows by the singularity analysis.
\section{Convex polyiamonds} \label{CP-section}
In this section, we enumerate the number of convex polyiamonds with respect to their perimeter. To that end, we first define a few notations.  We say that the bottom (resp. top) cell of the $i$-th column of a polyiamond $\nu$ is at the {\it position} $\ell$ if its lowest (resp. highest) point coincides with the line $y=\ell$. Let $b(\nu,i)$ and $u(\nu,i)$ denote the position of the bottom and top cells of the $i$-th column in the polyiamond $\nu$. Denote the number of columns of $\nu$ by $col(\nu)$. We define $\CP_k^{bu(i)}$ to be the set of all polyiamonds $\nu$ in $\CP$ where ($1$) the first column of $\nu$ is of type $i$ and has $k$ up-cells, and ($2$) $b(\nu,j+1)>b(\nu,j)$ and $u(\nu,j+1)<u(\nu,j)$ for all $j=1,2,\ldots,col(\nu)-1$. Also, we set $\CP_k^{u(i)}$ to be the set of polyiamonds where the condition ($2$) is replaced by the condition ($2'$) $u(\nu,j+1)<u(\nu,j)$ for all $j=1,2,\ldots,col(\nu)-1$. Similarly, we define $\CP_k^{b(i)}$ so that, ($2$) is replaced by the condition; ($2''$) $b(\nu,j+1)>b(\nu,j)$ for all $j=1,2,\ldots,col(\nu)-1$. Finally, let $\CP_k^{(i)}$ to be the set of convex polyiamonds $\nu$ satisfying only ($1$).

The enumeration of the convex polyiamonds with respect to perimeter is conducted as follows. The first step is to perform the counting in $\CP^{bu}:=\cup_{i=1}^4 \cup_{k\geq 1} \CP_k^{bu(i)}$ in a relatively straightforward manner. Then, we extend the result from $\CP^{bu}$ to $\CP^{u}:=\cup_{i=1}^4 \cup_{k\geq 1} \CP_k^{u(i)}$ and $\CP^{b}:=\cup_{i=1}^4 \cup_{k\geq 1} \CP_k^{b(i)}$. This can be done since the generating functions over these subsets can be written recursively in terms of the corresponding generating functions in $\CP^{bu}$. Finally, for the last step, we derive the results for $\CP$ by lifting up the result obtained for $\CP^u$, $\CP^b$, and $\CP^{bu}$.

\subsection{Enumeration over $\CP^{bu}$}
Suppose $F_k^{bu(i)}(p)$ is the generating function of polyiamonds in $\CP_k^{bu(i)}$ with respect to their perimeter; that is
\beq
F_k^{bu(i)}(p) := \sum_{\nu\in \CP_k^{bu(i)}} p^{per(\nu)}.
\feq
Moreover, define $F^{bu(i)}(p,u)=\sum_{k\geq0}F^{bu(i)}_k(p)u^k$ for all $i=1,2,3,4$, where $u$ marks the number of up-cells in the first column.

\begin{figure}[htp]
\begin{picture}(40,37)
\setlength{\unitlength}{.5mm} \linethickness{0.3mm}
\put(-20,67){Case (1)} \put(30,67){Case (2)}\put(80,67){Case (3)}
\def\uptt{\put(0,0){\line(1,0){10}}\put(10,0){\line(-2,3){5}}\put(0,0){\line(2,3){5}}}
\def\dowtt{\put(0,10){\line(1,0){10}}\put(10,10){\line(-2,-3){5}}\put(0,10){\line(2,-3){5}}}
\put(-40,10){\put(0,0){\uptt\put(5,-2.5){\dowtt}}
\put(5,7.5){\uptt\put(5,-2.5){\dowtt}}
\multiput(16,17)(1.6,2){3}{\circle*{1}}
\put(15,22.5){\uptt\put(5,-2.5){\dowtt}}
\put(20,30){\uptt\put(5,-2.5){\dowtt}}
\put(25,37.5){\uptt\put(5,-2.5){\dowtt}}
\put(10,15){\line(2,3){7}}\put(20,15){\line(2,3){7}}}
\put(-5,0){\put(0,10){\put(0,0){\uptt\put(5,-2.5){\dowtt}}
\put(5,7.5){\uptt\put(5,-2.5){\dowtt}}
\multiput(16,17)(1.6,2){3}{\circle*{1}}
\put(15,22.5){\uptt\put(5,-2.5){\dowtt}}
\put(20,30){\uptt\put(5,-2.5){\dowtt}}
\put(25,37.5){\uptt\put(5,-2.5){\dowtt}}
\put(10,15){\line(2,3){7}}\put(20,15){\line(2,3){7}}}
\put(15,17.5){\put(5,7.5){\uptt\put(5,-2.5){\dowtt}}
\multiput(16,17)(1.6,2){3}{\circle*{1}}
\put(15,22.5){\uptt\put(5,-2.5){\dowtt}}
\put(10,15){\line(2,3){7}}\put(20,15){\line(2,3){7}}
\multiput(25,37.5)(4,0){6}{\line(1,0){2}}
\multiput(-5,-7.5)(2,0){12}{\circle{0.1}}}}
\put(45,0){\put(0,10){\put(0,0){\uptt\put(5,-2.5){\dowtt}}
\put(5,7.5){\uptt\put(5,-2.5){\dowtt}}
\multiput(16,17)(1.6,2){3}{\circle*{1}}
\put(15,22.5){\uptt\put(5,-2.5){\dowtt}}
\put(20,30){\uptt\put(5,-2.5){\dowtt}}
\put(25,37.5){\uptt\put(5,-2.5){\dowtt}}
\put(10,15){\line(2,3){7}}\put(20,15){\line(2,3){7}}}
\put(15,17.5){\put(5,7.5){\uptt\put(5,-2.5){\dowtt}}
\multiput(16,17)(1.6,2){3}{\circle*{1}} \put(15,22.5){\uptt}
\put(10,15){\line(2,3){5}}\put(20,15){\line(2,3){5}}
\multiput(25,37.5)(4,0){6}{\line(1,0){2}}
\multiput(-5,-7.5)(2,0){12}{\circle{0.1}}}}
\end{picture}
\caption{Decomposition of a polyiamond in $\CP_k^{bu(1)}$ with
respect the first and second columns relative position} \label{figcpbu1}
\end{figure}
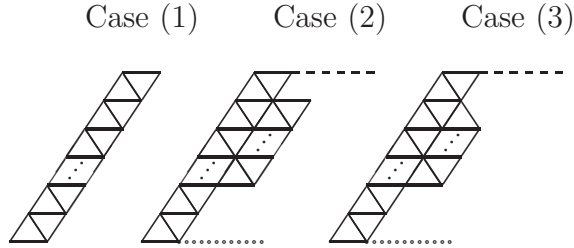

Note that any $\nu\in \CP_k^{bu(1)}$, depending on the type of the first and the second columns (if any) and their relative position, falls into one of cases described in Figure \ref{figcpbu1}.

Adding up the contributions of these terms, we arrive at
\beqn
F^{bu(1)}_k(p)&=p^{2k+2}+\sum_{j=1}^k(k+1-j)p^{2k+2-2j}(F^{bu(1)}_j(p)+F^{bu(2)}_{j-1}(p)),\label{Fbu-aa1}
\feqn
for all $k\geq1$.
In addition, by a cell growing/pruning based argument, similar to the one used in the previous sections, we obtain
\beqn \label{Fbu-lab}
&& F^{bu(2)}_k(p)=pF^{bu(1)}_k(p),\quad F^{bu(3)}_k(p)=\frac{1}{p}F^{bu(1)}_k(p),\quad k\geq1, \notag\\
&& F^{bu(4)}_k(p)=F^{bu(1)}_k(p) \quad F^{bu(2)}_0(p)=p^3, \quad F^{bu(i)}_0(p)=0, \ k\geq1,\,i=1,3,4.
\feqn
Multiplying by $u^k$ and adding up over $k\geq1$, the recursions \eqref{Fbu-aa1} and \eqref{Fbu-lab} can be written as the system of equations
\beqn \label{Fbu-sys}
&&F^{bu(1)}(p,u)=\frac{p^4u}{(1-p^2u)^2}+\frac{p^2}{(1-p^2u)^2}(F^{bu(1)}(p,u)+uF^{bu(2)}(p,u)),\notag\\
&&F^{bu(2)}(p,u)=p^3+pF^{bu(1)}(p,u), \notag\\
&& F^{bu(3)}(p,u)=\frac{1}{p}F^{bu(1)}(p,u),\quad  F^{bu(4)}(p,u)=F^{bu(1)}(p,u).
\feqn
By solving this system, we obtain the following result.
\begin{lemma}\label{lemcp1}
$F^{bu(1)}(p,u)=\frac{p^4u(1+p-p^2u)}{1-p^2-2p^2u-p^3u+p^4u^2}.$
\end{lemma}

\subsection{Enumeration over $\CP^u$}
Suppose $F_k^{u(i)}(p)$ is the generating function of polyiamonds in $\CP_k^{u(i)}$ with respect to their perimeter; that is
\beq
F_k^{u(i)}(p) := \sum_{\nu\in \CP_k^{u(i)}} p^{per(\nu)}.
\feq
Moreover, define $F^{u(i)}(p,u)=\sum_{k\geq0}F^{u(i)}_k(p)u^k$ for all $i=1,2,3,4$, where $u$ marks the number of up-cells in the first column.

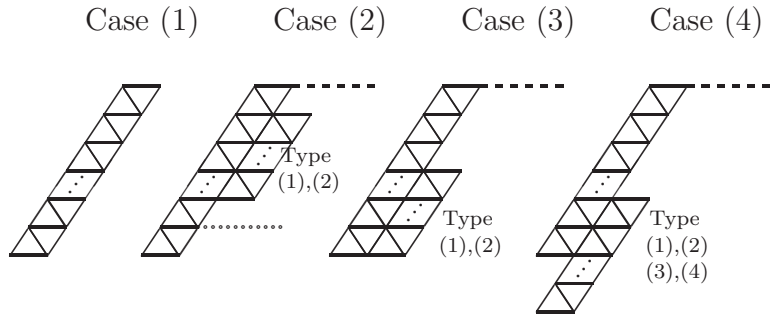
\begin{figure}[htp]
\begin{picture}(70,45)
\setlength{\unitlength}{.5mm} \linethickness{0.3mm}
\put(0,12){
\put(-20,70){Case (1)} \put(30,70){Case (2)}\put(80,70){Case (3)}\put(130,70){Case (4)}
\def\uptt{\put(0,0){\line(1,0){10}}\put(10,0){\line(-2,3){5}}\put(0,0){\line(2,3){5}}}
\def\dowtt{\put(0,10){\line(1,0){10}}\put(10,10){\line(-2,-3){5}}\put(0,10){\line(2,-3){5}}}
\put(-40,10){\put(0,0){\uptt\put(5,-2.5){\dowtt}}
\put(5,7.5){\uptt\put(5,-2.5){\dowtt}}
\multiput(16,17)(1.6,2){3}{\circle*{1}}
\put(15,22.5){\uptt\put(5,-2.5){\dowtt}}
\put(20,30){\uptt\put(5,-2.5){\dowtt}}
\put(25,37.5){\uptt\put(5,-2.5){\dowtt}}
\put(10,15){\line(2,3){7}}\put(20,15){\line(2,3){7}}}
\put(-5,0){\put(0,10){\put(0,0){\uptt\put(5,-2.5){\dowtt}}
\put(5,7.5){\uptt\put(5,-2.5){\dowtt}}
\multiput(16,17)(1.6,2){3}{\circle*{1}}
\put(15,22.5){\uptt\put(5,-2.5){\dowtt}}
\put(20,30){\uptt\put(5,-2.5){\dowtt}}
\put(25,37.5){\uptt\put(5,-2.5){\dowtt}}
\put(10,15){\line(2,3){7}}\put(20,15){\line(2,3){7}}}
\put(15,17.5){\put(5,7.5){\uptt\put(5,-2.5){\dowtt}}
\multiput(16,17)(1.6,2){3}{\circle*{1}}
\put(15,22.5){\uptt\put(5,-2.5){\dowtt}}
\put(10,15){\line(2,3){7}}\put(20,15){\line(2,3){7}}
\multiput(25,37.5)(4,0){6}{\line(1,0){2}}
\multiput(0,0)(2,0){12}{\circle{0.1}}}}
\put(32,35){\tiny Type}
\put(31,28){\tiny (1),(2)}
\put(45,0){\put(0,10){\put(0,0){\uptt\put(5,-2.5){\dowtt}}
\put(5,7.5){\uptt\put(5,-2.5){\dowtt}}
\multiput(16,17)(1.6,2){3}{\circle*{1}}
\put(15,22.5){\uptt\put(5,-2.5){\dowtt}}
\put(20,30){\uptt\put(5,-2.5){\dowtt}}
\put(25,37.5){\uptt\put(5,-2.5){\dowtt}}
\put(10,15){\line(2,3){7}}\put(20,15){\line(2,3){7}}}
\put(5,2.5){\put(5,7.5){\uptt\put(5,-2.5){\dowtt}}
\multiput(16,17)(1.6,2){3}{\circle*{1}} \put(15,22.5){\uptt\put(5,-2.5){\dowtt}}
\put(10,15){\line(2,3){5}}\put(20,15){\line(2,3){5}}
\multiput(35,52.5)(4,0){6}{\line(1,0){2}}}
\put(30,18){\tiny Type}
\put(29,11){\tiny (1),(2)}}
\put(100,0){\put(0,10){\put(0,0){\uptt\put(5,-2.5){\dowtt}}
\put(5,7.5){\uptt\put(5,-2.5){\dowtt}}
\multiput(16,17)(1.6,2){3}{\circle*{1}}
\put(15,22.5){\uptt\put(5,-2.5){\dowtt}}
\put(20,30){\uptt\put(5,-2.5){\dowtt}}
\put(25,37.5){\uptt\put(5,-2.5){\dowtt}}
\put(10,15){\line(2,3){7}}\put(20,15){\line(2,3){7}}}
\put(-5,-12.5){\put(5,7.5){\uptt\put(5,-2.5){\dowtt}}
\multiput(16,17)(1.6,2){3}{\circle*{1}} \put(15,22.5){\uptt\put(5,-2.5){\dowtt}}
\put(20,30){\uptt\put(5,-2.5){\dowtt}}
\put(10,15){\line(2,3){5}}\put(20,15){\line(2,3){5}}
\multiput(45,67.5)(4,0){6}{\line(1,0){2}}}
\put(30,18){\tiny Type}
\put(29,11){\tiny (1),(2)}
\put(29,4){\tiny (3),(4)}}
}
\end{picture}
\caption{Decomposition of a polyiamond in $\CP_k^{u(1)}$ with
respect to the relative position of the first and second columns. Only the
cases whose second column is of type $1$ are shown. The rest of cases are obtained by simply considering other types of second columns stated in the picture. }\label{figcpu1}
\end{figure}

Note that, inspecting the type of the first and second columns (if any) and their relative position, any $\nu\in \CP_k^{u(1)}$ belongs to one of cases described in Figure \ref{figcpu1}.

Adding up the contributions of these cases, we get
\noindent\beqn \label{Fuk1}
&&F^{u(1)}_k(p)=p^{2k+2}+\sum_{j=1}^{k-1}(k-j)p^{2k+2-2j}(F^{bu(1)}_j(p)+F^{bu(2)}_{j-1}(p))\notag\\ &&+\sum_{j=1}^kp^{2k+2-2j}(F^{u(1)}_j(p)+F^{u(2)}_{j-1}(p))\\
&&+\sum_{j=1}^k(p^{2k+2-2}+p^{2k+2-4}+\cdots+p^{2k+2-2(j-1)})(F^{u(1)}_j(p)+F^{u(2)}_{j-1}(p))\notag\\
&&+\sum_{j=1}^k(p^{2k+2-2}+p^{2k+2-4}+\cdots+p^{2k+2-2j})(F^{u(3)}_{j+1}(p)+F^{u(4)}_{j}(p))\notag\\
&&+\sum_{j\geq k+1}(p^{2k+2-2}+p^{2k+2-4}+\cdots+p^2)(F^{u(1)}_j(p) +F^{u(2)}_{j-1}(p)+F^{u(3)}_{j+1}(p)+F^{u(4)}_{j}(p)),\notag
\feqn
for $k\geq1$.
Once again, we apply the {\it cell growing/pruning based argument} and get
\beqn \label{Fuk234}
&&F^{u(2)}_k(p)=pF^{u(1)}_{k}(p),\quad F^{u(3)}_k(p)=\frac{1}{p}F^{u(1)}_{k}(p), \quad k\geq1,\notag \\
&&F^{u(4)}_k(p)=F^{u(1)}_{k}(p),\quad  F^{u(2)}_0(p)=p^3, \quad F^{u(i)}_0(p)=0, \quad k\geq1,\, i=1,3,4.
\feqn
Multiplying by $u^k$ and summing up over $k\geq1$, equations \eqref{Fuk1} and \eqref{Fuk234} can be written as
\beqn \label{Fu234qq}
&&F^{u(1)}(p,u)=\frac{p^4u}{1-p^2u}+\frac{p^4u}{(1-p^2u)^2}(F^{bu(1)}(p,u)+uF^{bu(2)}(p,u))\notag \\
&&\quad -\frac{p^2u}{(1-u)(1-p^2u)}(F^{u(1)}(p,u)+uF^{u(2)}(p,u)+\frac{1}{u}F^{u(3)}(p,u)+F^{u(4)}(p,u))\notag \\
&&\quad +\frac{p^2u}{(1-u)(1-p^2u)}(F^{u(1)}(p,1)+F^{u(2)}(p,1)+F^{u(3)}(p,1)+F^{u(4)}(p,1))
\feqn
Then by \eqref{Fbu-sys}, we get
\beqn
F^{u(1)}(p,u)&=\frac{p^4(1+p)u}{1-p^2u}+\frac{p^4u}{(1-p^2u)^2}(F^{bu(1)}(p,u)+uF^{bu(2)}(p,u))\nonumber \\
&-\frac{p}{(1-u)(1-p^2u)}((1+pu)^2F^{u(1)}(p,u)-u(1+p)^2F^{u(1)}(p,1)).\label{eqFu1qqqqq}
\feqn
Again, we apply the kernel method to \eqref{eqFu1qqqqq}. To that goal, let $u=u_0=\frac{1-p-\sqrt{1-2p-3p^2}}{2p^2}$ (which is the generating function for Motzkin numbers) and use Lemma \ref{lemcp1} and \eqref{Fbu-sys}. Then \eqref{eqFu1qqqqq} gives
\beqn \label{Fu1p1-sol}
F^{u(1)}(p,1)=\frac{p^3(3p^2+p-1+(p+1)\sqrt{1-2p-3p^2})}{(1-2p-4p^2)(1+p)}.
\feqn
By substituting \eqref{Fu1p1-sol} into \eqref{eqFu1qqqqq} and using Lemma \ref{lemcp1} and \eqref{Fbu-sys}, we obtain the following.
\begin{lemma}\label{lemFu1}
\begin{align*}
F^{u(1)}(p,u)&=\frac{p^4(p+1)(u\sqrt{1-2p-3p^2}+(1-p)u-2)}{(1-2p-4p^2)(p^2u^2+pu-u+1)}
-p^2\\
&+\frac{p^2((u+2)p^5u-(3u+2)p^4-(2u+4)p^3+(u+3)p^2+2p-1)}{(p^4u^2-p^3u-2p^2u-p^2+1)(4p^2+2p-1)}
\end{align*}
\end{lemma}
This lemma and \eqref{Fu234qq} imply that the generating function $\sum_{i=1}^4 F^{u(i)}(p,1)$ is:
\beq
\frac{p^2(-1+p+2p^2-p^3+(p+1)\sqrt{1-2p-3p^2})}{1-2p-4p^2}.
\feq

\begin{figure}[htp]
	\begin{picture}(60,45)
	\setlength{\unitlength}{.5mm} \linethickness{0.3mm}
	\put(0,0){
		\put(-20,78){Case (1)} \put(30,78){Case (2)}\put(80,78){Case (3)}\put(130,78){Case (4)}
		\def\uptt{\put(0,0){\line(1,0){10}}\put(10,0){\line(-2,3){5}}\put(0,0){\line(2,3){5}}}
		\def\dowtt{\put(0,10){\line(1,0){10}}\put(10,10){\line(-2,-3){5}}\put(0,10){\line(2,-3){5}}}
		\put(-40,10){\put(0,0){\uptt\put(5,-2.5){\dowtt}}
			\put(5,7.5){\uptt\put(5,-2.5){\dowtt}}
			\multiput(16,17)(1.6,2){3}{\circle*{1}}
			\put(15,22.5){\uptt\put(5,-2.5){\dowtt}}
			\put(20,30){\uptt\put(5,-2.5){\dowtt}}
			\put(25,37.5){\uptt\put(5,-2.5){\dowtt}}
			\put(10,15){\line(2,3){7}}\put(20,15){\line(2,3){7}}}
		\put(-5,0){\put(0,10){\put(0,0){\uptt\put(5,-2.5){\dowtt}}
				\put(5,7.5){\uptt\put(5,-2.5){\dowtt}}
				\multiput(16,17)(1.6,2){3}{\circle*{1}}
				\put(15,22.5){\uptt\put(5,-2.5){\dowtt}}
				\put(20,30){\uptt\put(5,-2.5){\dowtt}}
				\put(25,37.5){\uptt\put(5,-2.5){\dowtt}}
				\put(10,15){\line(2,3){7}}\put(20,15){\line(2,3){7}}}
			\put(10,10){\put(5,7.5){\uptt\put(5,-2.5){\dowtt}}
				\multiput(16,17)(1.6,2){3}{\circle*{1}}
				\put(15,22.5){\uptt\put(5,-2.5){\dowtt}}
				\put(10,15){\line(2,3){7}}\put(20,15){\line(2,3){7}}
				\multiput(25,37.5)(4,0){6}{\line(1,0){2}}
				\multiput(0,0)(2,0){12}{\circle{0.1}}}}
		\put(32,35){\tiny Type}
		\put(31,28){\tiny (1),(2)}
		\put(45,0){\put(0,10){\put(0,0){\uptt\put(5,-2.5){\dowtt}}
				\put(5,7.5){\uptt\put(5,-2.5){\dowtt}}
				\multiput(16,17)(1.6,2){3}{\circle*{1}}
				\put(15,22.5){\uptt\put(5,-2.5){\dowtt}}
				\put(20,30){\uptt\put(5,-2.5){\dowtt}}
				\put(25,37.5){\uptt\put(5,-2.5){\dowtt}}
				\put(10,15){\line(2,3){7}}\put(20,15){\line(2,3){7}}}
			\put(20,25){\put(5,7.5){\uptt\put(5,-2.5){\dowtt}}
				\multiput(16,17)(1.6,2){3}{\circle*{1}} \put(15,22.5){\uptt\put(5,-2.5){\dowtt}}
				\put(10,15){\line(2,3){5}}\put(20,15){\line(2,3){5}}
				\multiput(-10,-15.5)(2,0){12}{\circle{0.1}}}
			\put(40,38){\tiny Type}
			\put(39,31){\tiny (1),(2)}}
		\put(100,0){\put(0,10){\put(0,0){\uptt\put(5,-2.5){\dowtt}}
				\put(5,7.5){\uptt\put(5,-2.5){\dowtt}}
				\multiput(16,17)(1.6,2){3}{\circle*{1}}
				\put(15,22.5){\uptt\put(5,-2.5){\dowtt}}
				\put(20,30){\uptt\put(5,-2.5){\dowtt}}
				\put(25,37.5){\uptt\put(5,-2.5){\dowtt}}
				\put(10,15){\line(2,3){7}}\put(20,15){\line(2,3){7}}}
			\put(25,32.5){\put(5,7.5){\uptt\put(5,-2.5){\dowtt}}
				\multiput(16,17)(1.6,2){3}{\circle*{1}} \put(15,22.5){\uptt\put(5,-2.5){\dowtt}}
				\put(20,30){\uptt\put(5,-2.5){\dowtt}}
				\put(10,15){\line(2,3){5}}\put(20,15){\line(2,3){5}}
				\multiput(-15,-22.5)(2,0){12}{\circle{0.1}}}
			\put(47,48){\tiny Type}
			\put(46,41){\tiny (1),(2)}}
	}
	\end{picture}
	\caption{Decomposition of a polyiamond in $\CP_k^{b(1)}$ with
		respect to the relative position of the first and second columns. Only the
		cases whose second column is of type $1$ are shown. The rest of cases are obtained by simply considering other types of second columns stated in the picture.}\label{figcpb1}
\end{figure}
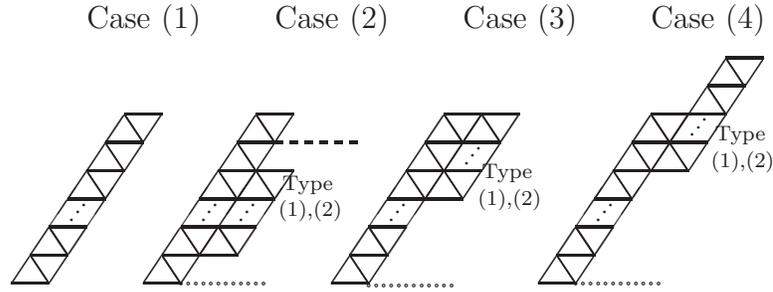

\subsection{Enumeration over $\CP^{b}$}
Suppose $F_k^{b(i)}(p)$ is the generating function of polyiamonds in $\CP_k^{b(i)}$ with respect to their perimeter; that is
\beq
F_k^{b(i)}(p) := \sum_{\nu\in \CP_k^{b(i)}} p^{per(\nu)}.
\feq
Moreover, define $F^{b(i)}(p,u)=\sum_{k\geq0}F^{b(i)}_k(p)u^k$ for all $i=1,2,3,4$, where $u$ marks the number of up-cells in the first column.

Note that any $\nu\in \CP_k^{b(1)}$, depending on the type of the first and second columns (if any) and their relative position, falls into one of cases described in Figure \ref{figcpb1}. Adding up the contributions of these terms, we get
\beq
&&F^{b(1)}_k(p)=p^{2k+2}+\sum_{j=1}^{k-1}(k-j)p^{2k+2-2j}(F^{bu(1)}_j(p)+F^{bu(2)}_{j-1}(p))\\ &&\quad+\sum_{j=1}^kp^{2k+2-2j}(F^{b(1)}_j(p)+F^{bu(2)}_{j-1}(p))\\
&&\quad +\sum_{j=1}^k(p^{2k+2-2}+p^{2k+3-4}+\cdots+p^{2k+2-2(j-1)})(F^{b(1)}_j(p)+F^{bu(2)}_{j-1}(p))\\
&&\quad +\sum_{j\geq k+1}(p^{2k+2-2}+p^{2k+2-4}+\cdots+p^2)(F^{b(1)}_j(p)+F^{bu(2)}_{j-1}(p)),
\feq
for all $k\geq1$.
An application of the {\it cell growing/pruning based argument} gives
\beq
&&F^{b(2)}_k(p)=F^{bu(2)}_{k}(p),\quad F^{b(3)}_k(p)=\frac{1}{p}F^{b(1)}_{k}(p),\quad k\geq1,\\
&&F^{b(4)}_k(p)=F^{bu(4)}_{k}(p),\quad F^{b(2)}_0(p)=p^3,\quad F^{b(i)}_0(p)=0,\quad k\geq1,\,  i=1,3,4.
\feq
By multiplying by $u^k$ and summing up over $k\geq1$, all these recursions can be written as
\beqn \label{qwqw}
&&F^{b(1)}(p,u)=\frac{p^4u}{1-p^2u}+\frac{p^4u}{(1-p^2u)^2}(F^{bu(1)}(p,u)+uF^{bu(2)}(p,u))\notag\\
&&+\frac{p^2u}{(1-u)(1-p^2u)}(F^{b(1)}(p,1)+F^{bu(2)}(p,1)-F^{b(1)}(p,u)-uF^{bu(2)}(p,u))\\
&&F^{b(2)}(p,u)=F^{bu(2)}(p,u),\  F^{b(3)}(p,u)=\frac{1}{p}F^{b(1)}(p,u), \ F^{b(4)}(p,u)=F^{bu(4)}(p,u). \notag
\feqn
This implies
\beqn
&&F^{b(1)}(p,u)=\frac{p^4u}{1-p^2u}+\frac{p^4u}{(1-p^2u)^2}(F^{bu(1)}(p,u)+uF^{bu(2)}(p,u))\notag \\
&&\quad +\frac{p^2u}{(1-u)(1-p^2u)}(F^{b(1)}(p,1)+F^{bu(2)}(p,1)-F^{b(1)}(p,u)-uF^{bu(2)}(p,u)).\label{eqFb1}
\feqn
We let $u=u_0=\frac{1-\sqrt{1-4p^2}}{2p^2}$ (the generating function for the Catalan numbers) in \eqref{eqFb1} and use Lemma \ref{lemcp1} and equation \eqref{Fbu-sys} to get
\beqn \label{eq123321}
&&F^{b(1)}(p,1)\\
&&=\frac{p^3((2p^2-1)(p^5+6p^3+5p^2-p-1)-(2p+1)(p^3-2p^2-p+1)(p+1)^2\sqrt{1-4p^2})}
{2(p+1)(5p^3+4p^2-p-1)(p^3-2p^2-p+1)}.\notag
\feqn
By substituting \eqref{eq123321} into \eqref{eqFb1} and using Lemma \ref{lemcp1}, we get
\begin{lemma}\label{lemFb1}
\beq
&& F^{b(1)}(p,u)=\frac{-u(p+1)(2p+1)p^5}{2(p^2u^2-u+1)(5p^3+4p^2-p-1)}\sqrt{1-4p^2}-p^2\\
&& +\frac{p^2(up^8+2(u+1)p^7+(u+3)p^6-(5u+6)p^5-(4u+7)p^4+(u+5)(1+p)p^2-p-1)}{(p^4u^2-p^3u-2p^2u-p^2+1)(5p^3+4p^2-p-1)}\\
&&-\frac{(u-2)(2p+1)(p+1)p^5}{(p^2u^2-u+1)(10p^3+8p^2-2p-2)},
\feq
\end{lemma}
An implication of this Lemma and \eqref{qwqw} is that
\beq
\sum_{i=1}^4 F^{b(i)}(p,1) =  \frac{p^2(2p^4-6p^3-7p^2+1-(2p+1)(p+1)^2\sqrt{1-4p^2})}{2(5p^3+4p^2-p-1)}.
\feq

\subsection{Enumeration over $\CP$}
Suppose $F_k^{(i)}(p)$ is the generating function of polyiamonds in $\CP_k^{(i)}$ with respect to their perimeter; that is
\beq
F_k^{(i)}(p) := \sum_{\nu\in \CP_k^{(i)}} p^{per(\nu)}.
\feq
Moreover, define $F^{(i)}(p,u)=\sum_{k\geq0}F^{(i)}_k(p)u^k$ for all $i=1,2,3,4$, where $u$ marks the number of up-cells in the first column.

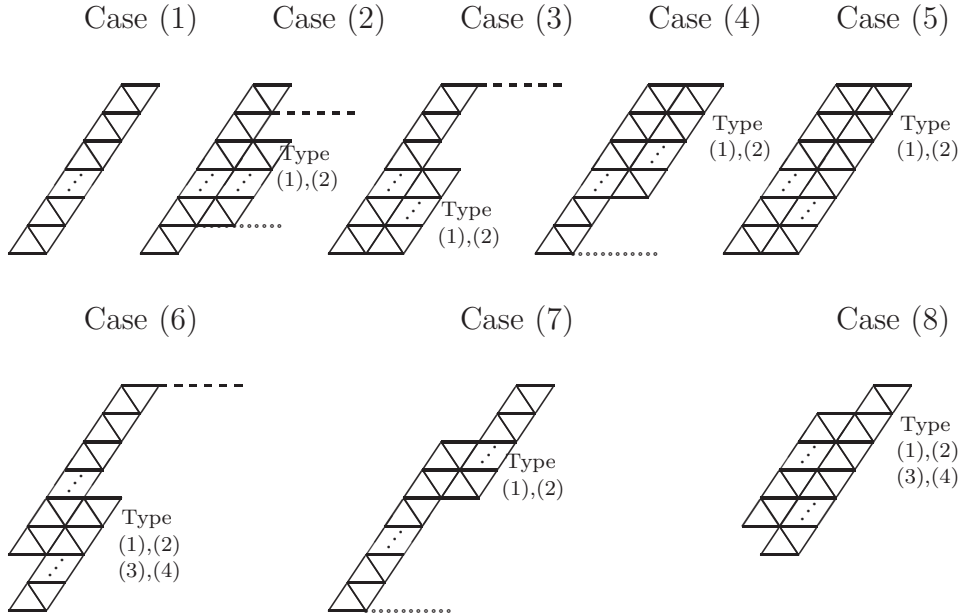
\begin{figure}[htp]
\begin{picture}(90,77)
\setlength{\unitlength}{.5mm} \linethickness{0.3mm}
\def\uptt{\put(0,0){\line(1,0){10}}\put(10,0){\line(-2,3){5}}\put(0,0){\line(2,3){5}}}
\def\dowtt{\put(0,10){\line(1,0){10}}\put(10,10){\line(-2,-3){5}}\put(0,10){\line(2,-3){5}}}
\put(0,80){
\put(-20,70){Case (1)} \put(30,70){Case (2)}\put(80,70){Case (3)}\put(130,70){Case (4)}\put(180,70){Case (5)}
\put(-40,10){\put(0,0){\uptt\put(5,-2.5){\dowtt}}
\put(5,7.5){\uptt\put(5,-2.5){\dowtt}}
\multiput(16,17)(1.6,2){3}{\circle*{1}}
\put(15,22.5){\uptt\put(5,-2.5){\dowtt}}
\put(20,30){\uptt\put(5,-2.5){\dowtt}}
\put(25,37.5){\uptt\put(5,-2.5){\dowtt}}
\put(10,15){\line(2,3){7}}\put(20,15){\line(2,3){7}}}
\put(-5,0){\put(0,10){\put(0,0){\uptt\put(5,-2.5){\dowtt}}
\put(5,7.5){\uptt\put(5,-2.5){\dowtt}}
\multiput(16,17)(1.6,2){3}{\circle*{1}}
\put(15,22.5){\uptt\put(5,-2.5){\dowtt}}
\put(20,30){\uptt\put(5,-2.5){\dowtt}}
\put(25,37.5){\uptt\put(5,-2.5){\dowtt}}
\put(10,15){\line(2,3){7}}\put(20,15){\line(2,3){7}}}
\put(10,10){\put(5,7.5){\uptt\put(5,-2.5){\dowtt}}
\multiput(16,17)(1.6,2){3}{\circle*{1}}
\put(15,22.5){\uptt\put(5,-2.5){\dowtt}}
\put(10,15){\line(2,3){7}}\put(20,15){\line(2,3){7}}
\multiput(25,37.5)(4,0){6}{\line(1,0){2}}
\multiput(5,7.5)(2,0){12}{\circle{0.1}}}}
\put(32,35){\tiny Type}
\put(31,28){\tiny (1),(2)}
\put(45,0){\put(0,10){\put(0,0){\uptt\put(5,-2.5){\dowtt}}
\put(5,7.5){\uptt\put(5,-2.5){\dowtt}}
\multiput(16,17)(1.6,2){3}{\circle*{1}}
\put(15,22.5){\uptt\put(5,-2.5){\dowtt}}
\put(20,30){\uptt\put(5,-2.5){\dowtt}}
\put(25,37.5){\uptt\put(5,-2.5){\dowtt}}
\put(10,15){\line(2,3){7}}\put(20,15){\line(2,3){7}}}
\put(5,2.5){\put(5,7.5){\uptt\put(5,-2.5){\dowtt}}
\multiput(16,17)(1.6,2){3}{\circle*{1}} \put(15,22.5){\uptt\put(5,-2.5){\dowtt}}
\put(10,15){\line(2,3){5}}\put(20,15){\line(2,3){5}}
\multiput(35,52.5)(4,0){6}{\line(1,0){2}}}
\put(30,20){\tiny Type}
\put(29,13){\tiny (1),(2)}}
\put(100,0){\put(0,10){\put(0,0){\uptt\put(5,-2.5){\dowtt}}
\put(5,7.5){\uptt\put(5,-2.5){\dowtt}}
\multiput(16,17)(1.6,2){3}{\circle*{1}}
\put(15,22.5){\uptt\put(5,-2.5){\dowtt}}
\put(20,30){\uptt\put(5,-2.5){\dowtt}}
\put(25,37.5){\uptt\put(5,-2.5){\dowtt}}
\put(10,15){\line(2,3){7}}\put(20,15){\line(2,3){7}}}
\put(15,17.5){\put(5,7.5){\uptt\put(5,-2.5){\dowtt}}
\multiput(16,17)(1.6,2){3}{\circle*{1}} \put(15,22.5){\uptt\put(5,-2.5){\dowtt}}
\put(20,30){\uptt\put(5,-2.5){\dowtt}}
\put(10,15){\line(2,3){5}}\put(20,15){\line(2,3){5}}
\multiput(-5,-7.5)(2,0){12}{\circle{0.1}}}
\put(47,43){\tiny Type}
\put(46,36){\tiny (1),(2)}}
\put(150,0){\put(0,10){\put(0,0){\uptt\put(5,-2.5){\dowtt}}
\put(5,7.5){\uptt\put(5,-2.5){\dowtt}}
\multiput(16,17)(1.6,2){3}{\circle*{1}}
\put(15,22.5){\uptt\put(5,-2.5){\dowtt}}
\put(20,30){\uptt\put(5,-2.5){\dowtt}}
\put(25,37.5){\uptt\put(5,-2.5){\dowtt}}
\put(10,15){\line(2,3){7}}\put(20,15){\line(2,3){7}}}
\put(5,2.5){\put(5,7.5){\uptt\put(5,-2.5){\dowtt}}
\multiput(16,17)(1.6,2){3}{\circle*{1}} \put(15,22.5){\uptt\put(5,-2.5){\dowtt}}
\put(20,30){\uptt\put(5,-2.5){\dowtt}}
\put(25,37.5){\uptt\put(5,-2.5){\dowtt}}
\put(30,45){\uptt\put(5,-2.5){\dowtt}}
\put(10,15){\line(2,3){5}}\put(20,15){\line(2,3){5}}}
\put(47,43){\tiny Type}
\put(46,36){\tiny (1),(2)}}
}
\put(0,0){
\put(-20,70){Case (6)} \put(80,70){Case (7)}\put(180,70){Case (8)}
\put(-40,0){\put(0,10){\put(0,0){\uptt\put(5,-2.5){\dowtt}}
\put(5,7.5){\uptt\put(5,-2.5){\dowtt}}
\multiput(16,17)(1.6,2){3}{\circle*{1}}
\put(15,22.5){\uptt\put(5,-2.5){\dowtt}}
\put(20,30){\uptt\put(5,-2.5){\dowtt}}
\put(25,37.5){\uptt\put(5,-2.5){\dowtt}}
\put(10,15){\line(2,3){7}}\put(20,15){\line(2,3){7}}}
\put(-5,-12.5){\put(5,7.5){\uptt\put(5,-2.5){\dowtt}}
\multiput(16,17)(1.6,2){3}{\circle*{1}} \put(15,22.5){\uptt\put(5,-2.5){\dowtt}}
\put(20,30){\uptt\put(5,-2.5){\dowtt}}
\put(10,15){\line(2,3){5}}\put(20,15){\line(2,3){5}}
\multiput(45,67.5)(4,0){6}{\line(1,0){2}}}
\put(30,18){\tiny Type}
\put(29,11){\tiny (1),(2)}
\put(29,4){\tiny (3),(4)}}
\put(45,-15){\put(0,10){\put(0,0){\uptt\put(5,-2.5){\dowtt}}
\put(5,7.5){\uptt\put(5,-2.5){\dowtt}}
\multiput(16,17)(1.6,2){3}{\circle*{1}}
\put(15,22.5){\uptt\put(5,-2.5){\dowtt}}
\put(20,30){\uptt\put(5,-2.5){\dowtt}}
\put(25,37.5){\uptt\put(5,-2.5){\dowtt}}
\put(10,15){\line(2,3){7}}\put(20,15){\line(2,3){7}}}
\put(25,32.5){\put(5,7.5){\uptt\put(5,-2.5){\dowtt}}
\multiput(16,17)(1.6,2){3}{\circle*{1}} \put(15,22.5){\uptt\put(5,-2.5){\dowtt}}
\put(20,30){\uptt\put(5,-2.5){\dowtt}}
\put(10,15){\line(2,3){5}}\put(20,15){\line(2,3){5}}
\multiput(-15,-22.5)(2,0){12}{\circle{0.1}}}
\put(47,48){\tiny Type}
\put(46,41){\tiny (1),(2)}}
\put(150,0){\put(5,17.5){\put(0,0){\uptt\put(5,-2.5){\dowtt}}
\put(5,7.5){\uptt\put(5,-2.5){\dowtt}}
\multiput(16,17)(1.6,2){3}{\circle*{1}}
\put(15,22.5){\uptt\put(5,-2.5){\dowtt}}
\put(10,15){\line(2,3){7}}\put(20,15){\line(2,3){7}}}
\put(5,2.5){\put(5,7.5){\uptt\put(5,-2.5){\dowtt}}
\multiput(16,17)(1.6,2){3}{\circle*{1}} \put(15,22.5){\uptt\put(5,-2.5){\dowtt}}
\put(20,30){\uptt\put(5,-2.5){\dowtt}}
\put(25,37.5){\uptt\put(5,-2.5){\dowtt}}
\put(30,45){\uptt\put(5,-2.5){\dowtt}}
\put(10,15){\line(2,3){5}}\put(20,15){\line(2,3){5}}}
\put(47,43){\tiny Type}
\put(46,36){\tiny (1),(2)}\put(46,29){\tiny (3),(4)}}
}
\end{picture}
\caption{Decomposition of a polyiamond in $\CP_k^{(1)}$ with respect to the relative position of the first and second columns. Only the cases whose second column is of type $1$ are shown. The rest of cases are obtained by simply considering other types of second columns stated in the picture.}\label{figcp1a}
\end{figure}

It is easy to see that each $\nu\in \CP_k^{(1)}$ belongs to one of cases described in Figure \ref{figcp1a}.

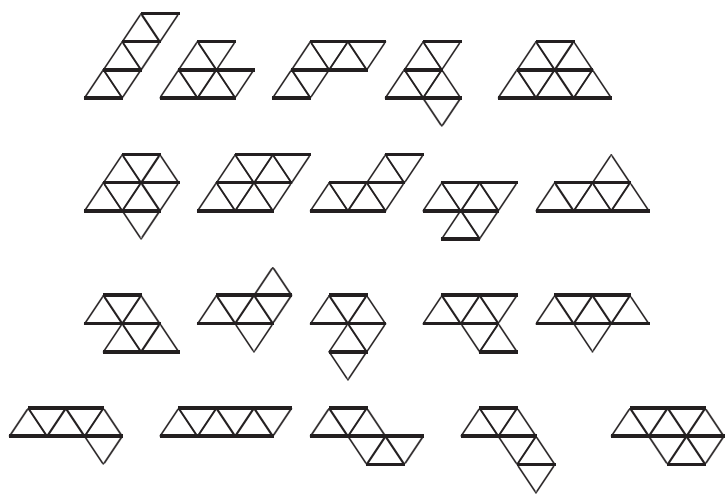
\begin{figure}[htp]
\begin{picture}(80,60)
\setlength{\unitlength}{.5mm} \linethickness{0.3mm}
\def\uptt{\put(0,0){\line(1,0){10}}\put(10,0){\line(-2,3){5}}\put(0,0){\line(2,3){5}}}
\def\dowtt{\put(0,10){\line(1,0){10}}\put(10,10){\line(-2,-3){5}}\put(0,10){\line(2,-3){5}}}
\put(0,100){
	\put(0,0){\uptt\put(5,-2.5){\dowtt}}\put(5,7.5){\uptt\put(5,-2.5){\dowtt}}\put(10,15){\uptt\put(5,-2.5){\dowtt}}
	\put(20,0){\put(0,0){\uptt\put(5,-2.5){\dowtt}}\put(5,7.5){\uptt\put(5,-2.5){\dowtt}}
		\put(10,0){\uptt\put(5,-2.5){\dowtt}}}
	\put(50,0){\put(0,0){\uptt\put(5,-2.5){\dowtt}}\put(5,7.5){\uptt\put(5,-2.5){\dowtt}}
		\put(15,7.5){\uptt\put(5,-2.5){\dowtt}}}
	\put(80,0){\put(0,0){\uptt\put(5,-2.5){\dowtt}}\put(5,7.5){\uptt\put(5,-2.5){\dowtt}}
		\put(10,0){\uptt\put(0,-10){\dowtt}}}
	\put(110,0){\put(0,0){\uptt\put(5,-2.5){\dowtt}}\put(5,7.5){\uptt\put(5,-2.5){\dowtt}}
		\put(10,0){\uptt\put(5,-2.5){\dowtt}\put(10,0){\uptt}\put(5,7.5){\uptt}}}
}
\put(0,70){
	\put(0,0){\put(0,0){\uptt\put(5,-2.5){\dowtt}}\put(5,7.5){\uptt\put(5,-2.5){\dowtt}}
		\put(10,0){\uptt\put(5,-2.5){\dowtt}\put(0,-10){\dowtt}\put(5,7.5){\uptt}}}
	\put(30,0){\put(0,0){\uptt\put(5,-2.5){\dowtt}}\put(5,7.5){\uptt\put(5,-2.5){\dowtt}}
		\put(10,0){\uptt\put(5,-2.5){\dowtt}\put(10,5){\dowtt}\put(5,7.5){\uptt}}}
	\put(60,0){\put(0,0){\uptt\put(5,-2.5){\dowtt}}
		\put(10,0){\uptt\put(5,-2.5){\dowtt}\put(10,5){\dowtt}\put(5,7.5){\uptt}}}
	\put(90,0){\put(0,0){\uptt\put(5,-2.5){\dowtt}}
		\put(5,-7.5){\uptt\put(5,-2.5){\dowtt}\put(10,5){\dowtt}\put(5,7.5){\uptt}}}
	\put(120,0){\put(0,0){\uptt\put(5,-2.5){\dowtt}}
		\put(10,0){\uptt\put(5,-2.5){\dowtt}\put(10,0){\uptt}\put(5,7.5){\uptt}}}
}
\put(0,40){
	\put(0,0){\put(0,0){\uptt\put(5,-2.5){\dowtt}}
		\put(5,-7.5){\uptt\put(5,-2.5){\dowtt}\put(10,0){\uptt}\put(5,7.5){\uptt}}}
	\put(30,0){\put(0,0){\uptt\put(5,-2.5){\dowtt}}
		\put(10,0){\uptt\put(5,-2.5){\dowtt}\put(0,-10){\dowtt}\put(5,7.5){\uptt}}}
	\put(60,0){\put(0,0){\uptt\put(5,-2.5){\dowtt}}
		\put(5,-7.5){\uptt\put(5,-2.5){\dowtt}\put(0,-10){\dowtt}\put(5,7.5){\uptt}}}
	\put(90,0){\put(0,0){\uptt\put(5,-2.5){\dowtt}}
		\put(5,-7.5){\put(5,-2.5){\dowtt}\put(10,5){\dowtt}\put(10,0){\uptt}\put(5,7.5){\uptt}}}
	\put(120,0){\put(0,0){\uptt\put(5,-2.5){\dowtt}}
		\put(5,-7.5){\put(5,-2.5){\dowtt}\put(10,5){\dowtt}\put(15,7.5){\uptt}\put(5,7.5){\uptt}}}
}
\put(-20,10){
	\put(0,0){\put(0,0){\uptt\put(5,-2.5){\dowtt}}
		\put(5,-7.5){\put(15,-2.5){\dowtt}\put(10,5){\dowtt}\put(15,7.5){\uptt}\put(5,7.5){\uptt}}}
	\put(40,0){\put(0,0){\uptt\put(5,-2.5){\dowtt}}
		\put(5,-7.5){\put(20,5){\dowtt}\put(10,5){\dowtt}\put(15,7.5){\uptt}\put(5,7.5){\uptt}}}
	\put(80,0){\put(0,0){\uptt\put(5,-2.5){\dowtt}}
		\put(5,-7.5){\put(5,7.5){\uptt}\put(5,-2.5){\dowtt}\put(10,0){\uptt}\put(15,-2.5){\dowtt}}}
	\put(120,0){\put(0,0){\uptt\put(5,-2.5){\dowtt}}
		\put(5,-7.5){\put(5,7.5){\uptt}\put(5,-2.5){\dowtt}\put(10,0){\uptt}\put(10,-10){\dowtt}}}
	\put(160,0){\put(0,0){\uptt\put(5,-2.5){\dowtt}\put(15,-2.5){\dowtt}\put(20,0){\uptt}}
		\put(5,-7.5){\put(5,7.5){\uptt}\put(5,-2.5){\dowtt}\put(10,0){\uptt}\put(15,-2.5){\dowtt}}}
}
\end{picture}
\caption{Convex polyiamonds in $\CP^{(1)}$ with perimeter 8}\label{figex20}
\end{figure}

Adding up the contributions of these terms, we get
\beq
&&F^{(1)}_k(p)=p^{2k+2}+p^2(F^{(1)}_{k}(p)+F^{(2)}_{k-1}(p))\\  &&+\sum_{j=1}^{k-2}(k-1-j)p^{2k+2-2j}(F^{bu(1)}_j(p)+F^{bu(2)}_{j-1}(p))\\
&&+\sum_{j=1}^{k-1}p^{2k+2-2j}(F^{u(1)}_j(p)+F^{b(1)}_j(p)+F^{u(2)}_{j-1}(p)+F^{bu(2)}_{j-1}(p))\\
&&+\sum_{j=1}^k(p^{2k+2-2}+p^{2k+2-4}+\cdots+p^{2k+2-2(j-1)})(F^{u(1)}_{j}(p)+F^{u(2)}_{j-1}(p))\\
&&+\sum_{j=1}^k(p^{2k+2-2}+p^{2k+2-4}+\cdots+p^{2k+2-2j})(F^{u(3)}_{j+1}(p)+F^{u(4)}_{j}(p))\\
&&+\sum_{j\geq k+1}(p^{2k+2-2}+p^{2k+2-4}+\cdots+p^2)(F^{u(1)}_{j}(p)+F^{u(2)}_{j-1}(p))\\
&&+\sum_{j\geq k+1}(p^{2k+2-2}+p^{2k+2-4}+\cdots+p^2)(F^{u(3)}_{j+1}(p)+F^{u(4)}_{j}(p))\\
&&+\sum_{j=1}^k(p^{2k+2-2}+p^{2k+2-4}+\cdots+p^{2k+2-2(j-1)})(F^{b(1)}_{j}(p)+F^{bu(2)}_{j-1}(p))\\
&&+\sum_{j\geq k+1}(p^{2k+2-2}+p^{2k+2-4}+\cdots+p^2)(F^{b(1)}_{j}(p)+F^{bu(2)}_{j-1}(p))\\
&&+\sum_{j\geq k+1}(j-1-k)p^2(F^{(1)}_{j}(p)+F^{u(2)}_{j-1}(p)) +\sum_{j\geq k+1}(j-k)p^2(F^{(3)}_{j+1}(p)+F^{u(4)}_{j}(p)),
\feq
for all $k\geq1$. As before, we apply the {\it cell growing/pruning based argument} to get
\beqn\label{Fikp-eqs}
&&F^{(2)}_k(p)=pF^{u(1)}_k(p), \quad F^{(3)}_k(p)=\frac{1}{p}F^{(1)}_k(p),\quad k\geq1,\notag\\
&&F^{(4)}_k(p)=F^{u(1)}_k(p),\quad F^{(2)}_0(p)=p^3, \quad F^{(i)}_0(p)=0,\quad  k\geq1,\,i=1,3,4.
\feqn
By multiplying equations in \eqref{Fikp-eqs} by $u^k$ and summing up over $k\geq1$, we obtain
\beqn \label{Fikp-eqs2}
F^{(2)}(p,u)=p^3+pF^{u(1)}(p,u),\quad
 F^{(3)}(p,u)=\frac{1}{p}F^{(1)}(p,u),\quad
F^{(4)}(p,u)=F^{u(1)}(p,u),
\feqn
and
\begin{align*}
&(1+p)(1-p-\frac{p}{(1-u)^2})F^{(1)}(p,u)=\frac{p^4u}{1-p^2u}+\frac{p^6u^2}{(1-p^2u)^2}(F^{bu(1)}(p,u)+uF^{bu(2)}(p,u))\\
&-\frac{(p^2u^2-p^2u+1)p^2}{(1-u)(1-p^2u)}(F^{u(1)}(p,u)+uF^{u(2)}(p,u)+F^{b(1)}(p,u)+uF^{bu(2)}(p,u))\\
&+p^2uF^{(2)}(p,u)-\frac{p^2u}{(1-u)(1-p^2u)}(\frac{1}{u}F^{u(3)}(p,u)+F^{u(4)}(p,u))\\
&+\frac{p^2u}{(1-u)(1-p^2u)}(F^{u(1)}(p,1)+F^{u(2)}(p,1)+F^{u(3)}(p,1)+F^{u(4)}(p,1))\\
&+\frac{p^2u}{(1-u)(1-p^2u)}(F^{b(1)}(p,1)+F^{bu(2)}(p,1))\\
&-\frac{p^2u(2-u)}{(1-u)^2}F^{(1)}(p,1)
+\frac{2p^2u}{1-u}\frac{\partial}{\partial u}F^{(1)}(p,u)\mid_{u=1}+\frac{p^2u}{1-u}\frac{\partial}{\partial u}F^{u(4)}(p,u)\mid_{u=1}\\
&+\frac{p^2u}{(1-u)^2}F^{(2)}(p,u)-\frac{p^2u(2-u)}{(1-u)^2}F^{(2)}(p,1)
+\frac{p^2u}{1-u}\frac{\partial}{\partial u}F^{(2)}(p,u)\mid_{u=1}\\ &+\frac{p^2u}{(1-u)^2}F^{u(4)}(p,u)
-\frac{pu}{(1-u)^2}F^{(1)}(p,1)-\frac{p^2u}{(1-u)^2}F^{u(4)}(p,1).
\end{align*}
We set $u^\pm=v_1=1\pm\sqrt{\frac{p}{1-p}}$ and use Lemmas \ref{lemcp1}, \ref{lemFu1} and \ref{lemFb1}, and obtain a system of two equations with two variables $F^{(1)}(p,1)$ and $\frac{\partial}{\partial u}F^{(1)}(p,u)\mid_{u=1}$. Solving this system, we get the following result. We omit the lengthy calculations for the sake of brevity.
\begin{lemma} \label{last-lem}
The generating function $F^{(1)}(p,1)$ is given by
\beq
&&\frac{p^4(-2p^5+p^4+3p^3-3p^2-2p+2)}{(4p^2+2p-1)(4p^5-7p^3+p^2+4p-1)}(p-\sqrt{1-2p-3p^2})\\
&&+\frac{p^6(p^2-1)(2p+1)}{2(2p^4-p^3-2p^2+p+1)(5p^3+4p^2-p-1)(2p-1)}\sqrt{1-4p^2}\\
&&-\frac{9}{5}p^2-\frac{13}{50}p-\frac{1577}{1000}+\frac{1}{p+1}\\
&&+\frac{570p^9+1009p^8-671p^7-1764p^6+210p^5+1251p^4-4p^3-459p^2-20p+63}
{88(4p^{10}-8p^9+p^8+19p^7-8p^6-20p^5+9p^4+10p^3-5p^2-2p+1)}\\
&&+\frac{24p^4+7p^3-p^2+8p-3)}{88(4p^5-7p^3+p^2+4p-1)}
+\frac{-p^2+19p+24}{500(5p^3+4p^2-p-1)}
+\frac{p^3+2p^2-p-1}{8(2p^4-p^3-2p^2+p+1)}\\
&&=p^4+p^5+3p^6+8p^7+{\bf20}p^8+58p^9+152p^{10}+427p^{11}+1155p^{12}+3211p^{13}+O(p^{14}).
\feq
\end{lemma}
Note that, for instance, the coefficient of $p^8$ is $20$, indicating that there are $20$  polyiamonds in $\CP^{(1)}$ of perimeter $8$ as seen in Figure \ref{figex20}.

By Lemma \ref{last-lem}, \eqref{Fu1p1-sol} and \eqref{Fikp-eqs2}, we obtain the first part of Theorem \ref{mthcp}.
Finally, when $p$ is near the dominant singularity of the generating function $F^{(1)}(p,1)$ (see Lemma \ref{last-lem}),  $F^{(1)}(p,1)$ can be approximated by
\beq
\frac{2114724967}{3292828497}-\frac{80}{15219}\sqrt{5}-\frac{2560}{441\sqrt{3}}\sqrt{1-3p}+O((1-3p)),
\feq
which, by the singularity analysis, gives the second part of Theorem \ref{mthcp}.
\medskip

\noindent{\bf Acknowledgements}: The authors would like to thank the referee for pointing out the mistake in the statement of Theorem 1.2 in the original draft.


\end{document}